\documentclass[11pt]{amsart}
\usepackage{amsxtra}
\usepackage{amssymb}
\theoremstyle{plain}

\newcommand{\cleqn}{\setcounter{equation}{0}}
\newcommand{\clth}{\setcounter{theorem}{0}}
\newcommand {\sectionnew}[1]{\section{#1}\cleqn\clth}

\newtheorem{theorem}{Theorem}[section]
\newtheorem{lemma}[theorem]{Lemma}
\newtheorem{definition-theorem}[theorem]{Definition-Theorem}
\newtheorem{proposition}[theorem]{Proposition}
\newtheorem{corollary}[theorem]{Corollary}
\newtheorem{definition}[theorem]{Definition}
\newtheorem{example}[theorem]{Example}
\newtheorem{remark}[theorem]{Remark}
\newtheorem{conjecture}[theorem]{Conjecture}
\newtheorem{notation}[theorem]{Notation}
\newcommand \bth[1] { \begin{theorem}\label{t#1} }
\newcommand \ble[1] { \begin{lemma}\label{l#1} }

\newcommand \bpr[1] { \begin{proposition}\label{p#1} }
\newcommand \bco[1] { \begin{corollary}\label{c#1} }
\newcommand \bde[1] { \begin{definition}\label{d#1}\rm }
\newcommand \bex[1] { \begin{example}\label{e#1}\rm }
\newcommand \bre[1] { \begin{remark}\label{r#1}\rm }
\newcommand \bcj[1] { \begin{conjecture}\label{j#1}\rm }

\newcommand \bnota[1] { \begin{notation}\label{n#1}\rm }
\renewcommand {\eth} { \end{theorem} }
\newcommand {\ele} { \end{lemma} }

\newcommand {\epr} { \end{proposition} }
\newcommand {\eco} { \end{corollary} }
\newcommand {\ede} { \end{definition} }
\newcommand {\eex} { \end{example} }
\newcommand {\ere} { \end{remark} }
\newcommand {\ecj} { \end{conjecture} }

\newcommand {\enota} { \end{notation} }
\newcommand \thref[1]{Theorem \ref{t#1}}
\newcommand \leref[1]{Lemma \ref{l#1}}
\newcommand \prref[1]{Proposition \ref{p#1}}
\newcommand \coref[1]{Corollary \ref{c#1}}

\newcommand \lb[1]{\label{#1}}

\def \Cset {{\mathbb C}}
\def \KK {{\mathbb K}}
\def \Zset {{\mathbb Z}}
\def \Tset {{\mathbb T}}
\def \Nset {{\mathbb N}}
\def \Qset {{\mathbb Q}}

\def \B  {{\mathcal{B}}}

\def \UU {{\mathcal{U}}}
\def \RR {{\mathcal{R}}}

\def \De {\Delta}   
\def \de {\delta}
\def \al {\alpha}
\def \be {\beta}
\def \la {\lambda}

\def \Om {\Omega} 
\def \om {\omega}
\def \ga {\gamma}
\def \de {\delta}
\def \Ga {\Gamma}
\def \Sig {\Sigma}

\def \mt  {\mapsto}

\def \hra {\hookrightarrow}

\def \ci  {\circ}

\def \rcor {\rangle}
\def \lcor {\langle}

\def \ol {\overline}

\def \id { {\mathrm{id}} }

\def \rank { {\mathrm{rank}} }

\def \g  {\mathfrak{g}}   

\def \sl {\mathfrak{sl}}

\def \n  {\mathfrak{n}}

\def \sl {\mathfrak{sl}}

\def \Mmn {M_{m,n}}

\DeclareMathOperator \GKdim { {\mathrm{GK \, dim}}}
\DeclareMathOperator \Span { {\mathrm{Span}} }
\DeclareMathOperator \Aut { {\mathrm{Aut}} }

\DeclareMathOperator \Wt { {\mathrm{wt}} }

\DeclareMathOperator \Char { {\mathrm{char}} }
\DeclareMathOperator \Ext  { {\mathrm{Ext}} }

\DeclareMathOperator \Hw  { {\mathrm{hw}} }

\newcommand \Spec { {\mathrm{Spec}} }
\begin{document}
\title[A conjecture of Goodearl and Lenagan]
{A proof of the Goodearl--Lenagan polynormality conjecture}
\author[Milen Yakimov]{Milen Yakimov}
\address{
Department of Mathematics \\
Louisiana State University \\
Baton Rouge, LA 70803 \\
U.S.A.
}
\email{yakimov@math.lsu.edu}
\date{}
\begin{abstract}
The quantum nilpotent algebras $\UU^w_-(\g)$, defined by
De Concini--Kac--Procesi and Lusztig, are 
large classes of iterated skew polynomial rings
with rich ring theoretic structure.
In this paper, we prove in an explicit way that all 
torus invariant prime ideals of the algebras $\UU^w_-(\g)$ are polynormal. In the 
special case of the algebras of quantum matrices, this construction yields  
explicit polynormal generating sets consisting of quantum minors 
for all of their torus invariant prime ideals. This gives a constructive
proof of the Goodearl--Lenagan polynormality conjecture \cite{GLen}.
Furthermore we prove that $\Spec \UU^w_-(\g)$ is normally separated for all
simple Lie algebras $\g$ and Weyl group elements $w$, 
and deduce from it that all algebras $\UU^w_-(\g)$ are catenary. 
\end{abstract}
\maketitle
\sectionnew{Introduction}
\lb{intro}
The algebras of quantum matrices $R_q[M_{m,n}]$ are algebras over a field 
$\KK$ generated by a rectangular array of generators 
$\{ x_{ij} \mid i=1, \ldots, m, j=1, \ldots, n \}$, such that 
rows and columns generate quantum affine space algebras. Generators
along antidiagonals commute, while generators along diagonals satisfy the 
commutation relation
\[
x_{ij} x_{lk} - x_{lk} x_{ij} = (q-q^{-1}) x_{ik} x_{lj}, \quad 
i < l, j<k,
\]
where $q \in \KK^* = \KK \backslash \{ 0 \}$. There is a canonical action 
of the $\KK$-torus $\Tset^{m+n} = (\KK^*)^{\times (m+n)}$ 
on $R_q[M_{m,n}]$ by algebra automorphisms. The classical minors have 
natural quantum analogs, which are elements of the algebras $R_q[M_{m,n}]$, 
see Section \ref{GL}.

The ring theoretic properties of the algebras of quantum matrices 
$R_q[M_{m,n}]$ have been heavily investigated since the mid 90's. In this paper 
we address a conjecture of Goodearl and Lenagan \cite{GLen} that all 
$\Tset^{m+n}$-invariant prime ideals of $R_q[M_{m,n}]$ have polynormal
generating sets consisting of quantum minors. This conjecture has 
been established \cite{GLen} only for $\min(m,n) \leq 2$ and $m=n=3$. Shortly after 
\cite{GLen}, Launois proved in \cite{La} that the 
$\Tset^{m+n}$-invariant prime ideals of $R_q[M_{m,n}]$ 
are generated by the quantum minors contained in those ideals 
when the base field has characteristic $0$ and 
$q$ is transcendental over $\Qset$. In \cite{Y} the author constructed 
explicit generating sets for these ideals under the same conditions on 
$\KK$ and $q$. These generating sets consist of quantum minors but do not contain all 
quantum minors in a given ideal. In \cite{Cas,GLL}
Casteels, Goodearl, Launois, and Lenagan determined 
what quantum minors belong to each $\Tset^{m+n}$-invariant 
prime ideal of $R_q[M_{m,n}]$, and completed the line of argument of 
\cite{La} for ideal generators (which results in larger generating 
sets than those in \cite{Y}).
In \cite{Y} we also constructed
explicit generating sets for torus invariant prime ideals of 
a much larger family of quantum nilpotent algebras. Meanwhile, there was 
no progress on the main part of the conjecture in \cite{GLen} on 
polynormal generating sequences. 

In this paper we give a constructive proof of the Goodearl--Lenagan
polynormality conjecture \cite{GLen} -- we construct explicit 
polynormal generating sets consisting of quantum minors  
for all $\Tset^{m+n}$-invariant 
prime ideals of $R_q[M_{m,n}]$. These sets are 
smaller than the sets of all quantum minors in a given $\Tset^{m+n}$-prime 
ideal. This is done for base fields of characteristic $0$ and when 
$q$ is transcendental over $\Qset$. In another theorem we
construct polynormal generating sets for all $\Tset^{m+n}$-invariant 
prime ideals of $R_q[M_{m,n}]$ under the weaker assumptions that $\KK$ 
is arbitrary and $q \in \KK^*$ is not a root of unity. In the second result 
the generating sets are bigger than those in the first one. The two results 
are proved by the same type of argument.

The algebras of quantum matrices are a subfamily of a much larger 
family of iterated skew polynomial rings defined 
by De Concini, Kac, and Procesi \cite{DKP}, and Lusztig \cite{L}.
These algebras $\UU^w_-$ are subalgebras of the negative part of 
the quantized universal enveloping algebra $\UU_q(\g)$ of an arbitrary 
simple Lie algebra $\g$ and are parametrized by elements $w$ 
of the Weyl group $W$ of $\g$. They are deformations of the universal
enveloping algebras $\UU( \n_- \cap w (\n_+))$, where $\n_\pm$ denote 
the nilradicals of a pair
of opposite Borel subalgebras of $\g$. 
These algebras can be considered as quantizations of the coordinate rings of 
Schubert cells. The geometry of the corresponding Poisson 
structures on all (partial) flag varieties was investigated 
by Brown, Goodearl, and the author in \cite{BGY,GY}. There is a 
canonical action of the $\KK$-torus $\Tset^r= (\KK^*)^{\times r}$ 
$(r = \rank \, \g)$ 
on $\UU_q(\g)$ by algebra automorphisms, which preserves all subalgebras $\UU^w_-$. 
The $\Tset^r$-invariant prime ideals of $\UU^w_-$ were
classified by M\`eriaux and Cauchon \cite{MC}, and the 
author \cite{Y}. Moreover, \cite{Y} described 
the poset structure of the $\Tset^r$-spectrum of 
the algebras $\UU^w_-$ and gave an explicit description 
of all of their torus invariant prime ideals. By the 
general spectrum stratification theorem \cite{GL} 
of Goodearl and Letzter, $\Spec \UU^w_-$ stratifies 
into a disjoint union of tori, which are parametrized 
by the $\Tset^r$-primes of $\UU^w_-$. Bell, Casteels, 
and Launois \cite{BCL} and the author \cite{Y3,Y5} computed 
the dimensions of these strata.

In this respect, it is interesting to understand what is the 
analog of the Goodearl--Lenagan polynormality conjecture for the 
$\Tset^r$-spectra of all algebras $\UU^w_-$. This is 
the next problem addressed in the paper. The
simple Lie algebras, which are not of $A$ type, have 
fundamental representations that are not minuscule. 
Thus one cannot expect to find generating sets 
for $\Tset^r$-primes of $\UU^w_-$, which consist of elements 
derived from the braid group orbits of highest weight
vectors of fundamental $\UU_q(\g)$-modules, since these sets 
will be too small. We prove that for an arbitrary base field
$\KK$ of characteristic $0$ and $q \in \KK^*$ which is
transcendental over $\Qset$, all $\Tset^r$-prime 
ideals of $\UU^w_-$ have explicit polynormal generating sets 
parametrized by certain weight vectors of the Demazure 
modules corresponding to fundamental weights and the given
Weyl group element $w$. In the special case when $w$ equals 
the longest element $w_0$ of the Weyl group $W$ and $\KK= \Cset(q)$ 
implicit polynormal generating sets 
were obtained by Caldero \cite{Cal}. 
In the $A$ case our generating sets 
consist of elements which are generalizations of quantum minors 
(even for arbitrary Weyl group elements $w$).
Furthermore, in the case of an arbitrary base field $\KK$ 
and $q \in \KK^*$ which is not a root of unity,
for all $\Tset^r$-prime ideals of $\UU^w_-$ we construct 
polynormal generating sets, which are parametrized by certain 
weight vectors in possibly higher Demazure modules. Therefore 
all $\Tset^r$-prime ideals of the algebras $\UU^w_-$ are 
polynormal for an arbitrary base field $\KK$ and $q \in \KK^*$ 
not a root of unity. In fact we prove stronger equivariant 
polynormality with respect to an action of the weight lattice
$P$ of $\g$, naturally embedded in $\Tset^r$. We refer 
to \S \ref{3.1} for the definitions of equivariant 
normality and polynormality, and 
to Theorems \ref{tpolynormal} and \ref{tpol2} for 
the precise statements. 
Our proofs 
of these results and the results for quantum matrices rely on 
theorems of Joseph \cite{J1,J2} and Gorelik \cite{G}. These facts 
appear in Sections \ref{hsp} and \ref{GL}. Brown and Goodearl 
constructed in \cite{BG0} polynormal generating sets for all 
torus invariant prime ideals of the quantum algebras 
of functions on simple groups. The difference with our 
situation is that the definition of the algebras 
$\UU^w_-$ is not in terms of quantum function algebras and
one cannot use $R$-matrix type commutation relations. 
As a result, no relation between polynormality for the 
algebras of quantum matrices and the quantum function 
algebras on simple groups was previously observed. 
Only after one realizes the algebras $\UU^w_-$ 
as quantum coordinate rings, this approach to polynormality 
becomes possible. This realization of the algebras $\UU^w_-$ coincides
up to an antiisomorphism
with Joseph's algebras $S^w_+$ \cite[\S 10.3.1]{J} (see \S 2.4),
which played an important role in his study \cite{J1,J} of the spectra
of quantum groups.

In Section \ref{catenarity} we prove that the spectra of all algebras 
$\UU^w_-$ are normally separated for an 
arbitrary base field $\KK$ and $q \in \KK^*$ not a root of unity.
The special case of the algebras of quantum matrices 
is due to Cauchon \cite{Cau}. The case when $w= w_0$ and $\KK= \Cset(q)$
was obtained by Caldero \cite{Cal}.
We give two proofs of this result. The first applies directly 
the polynormality and the second relies on a result of Gorelik 
\cite{G}. For both we use a theorem of Goodearl \cite{G1} 
to pass from graded normal separation of the $\Tset^r$-spectrum
of $\UU^w_-$ to normal separation of $\Spec \UU^w_-$. 

A celebrated theorem of Gabber establishes that the universal 
enveloping algebras of all solvable Lie algebras are catenary. 
Since the algebras $\UU^w_-$ are quantized universal enveloping 
algebras of nilpotent Lie algebras, one can conjecture that all 
of them should be catenary. We prove this here. The special case 
of catenarity of the algebras of quantum matrices was previously 
proved by Cauchon \cite{Cau}. The case when $w=w_0$ and $\Char \KK =0$ 
is due to Malliavin \cite{M}, Goodearl and Lenagan \cite{GLen0}.
We use a general result of Goodearl and Lenagan \cite{GLen0}, motivated 
by Gabber's work, 
that establishes that catenarity follows from normal separation 
and certain homological conditions, and a result 
of Levasseur and Stafford \cite{LS} that proves the latter for skew 
polynomial extensions. Finally, we derive explicit formulas for 
the heights of all $\Tset^r$-invariant prime ideals of the algebras 
$\UU^w_-$. This is done in Section \ref{catenarity}.

The algebras $\UU^w_-$ belong to the large class 
of so called (torsion free) Cauchon--Goodearl--Letzter 
(CGL) extensions 
\cite[Definition 3.1]{LLR}. 
The latter are iterated skew polynomial rings with 
a compatible torus action for which the Goodearl--Letzter 
stratification result produces a finite stratification and 
the Cauchon method of deleted derivations applies. It is 
very interesting, yet very difficult to prove or disprove 
whether all (torsion free) CGL extensions have normal 
separation and thus are catenary. One should note that 
general skew polynomial rings do not have this 
property as shown by Bell and Sigurdsson 
\cite[Example 2.10]{BS}.

We finish with a notational convention. All algebras $A$
which we consider in this paper are noetherian over 
an infinite field $\KK$. When such algebras are equipped 
with a rational $\Tset^r$-action 
by algebra automorphisms, a result of Brown and Goodearl 
\cite[Proposition II.2.9]{BG} applies to give that
all of their $\Tset^r$-primes are $\Tset^r$-invariant 
primes. (We refer to \cite[Sect. II.2]{BG} for a discussion of rational 
torus actions over arbitrary fields.) Because 
of this fact, we will use the two terms interchangeably.
The corresponding $\Tset^r$-spectrum will be denoted 
by $\Tset^r-\Spec A$.
\\ \hfill \\
{\bf Acknowledgements.} I am grateful to Ken Goodearl and 
St\'ephane Launois for their helpful comments and for sharing their 
knowledge of the existing literature. I would also like to thank 
Ken Goodearl for communicating his proof \cite{G2}
of \prref{he} to me and the referee whose valuable 
comments and suggestions helped me to improve the exposition. 

The research of the author was supported by National 
Science Foundation grants DMS-0701107 and DMS-1001632.
\sectionnew{Quantized nilpotent algebras}
\lb{qalg}
\subsection{}
\label{2.1}
In this section we set up our notation for quantized universal 
enveloping algebras and quantum function algebras, and recall 
past result on their spectra which will be needed
in the paper. 

Fix an arbitrary base field $\KK$ and $q \in \KK^*$ which is not a root 
of unity. Let $\g$ be a simple Lie algebra of rank $r$ with Cartan matrix
$(c_{ij})$. Denote by 
$\UU_q(\g)$ the quantized universal enveloping algebra of $\g$ over the 
base field $\KK$ with deformation parameter $q$. It is a Hopf algebra over 
$\KK$ with generators
\[
X^\pm_i, K_i^{\pm 1}, \; i=1, \ldots, r
\]
and relations
\begin{gather*}
K_i^{-1} K_i = K_i K^{-1}_i = 1, \, K_i K_j = K_j K_i, 
\\
K_i X^\pm_j K^{-1}_i = q_i^{\pm c_{ij}} X^\pm_j,
\\
X^+_i X^-_j - X^-_j X^+_i = \de_{i,j} \frac{K_i - K^{-1}_i}
{q_i - q^{-1}_i},
\\
\sum_{k=0}^{1-c_{ij}}
(-1)^k
\begin{bmatrix} 
1-c_{ij} \\ k
\end{bmatrix}_{q_{i}}
      (X^\pm_i)^k X^\pm_j (X^\pm_i)^{1-c_{ij}-k} = 0, \, i \neq j.
\end{gather*}
Here $q_i = q^{d_i}$ and $\{ d_1, \ldots, d_r \}$ is
the vector of positive relatively 
prime integers for which the matrix $(d_i c_{ij})$ is symmetric.  
The comultiplication, antipode, and counit of $\UU_q(\g)$ are given 
by:
\begin{multline*}
\De(K_i) = K_i \otimes K_i, \, \, 
\De(X^+_i) = X^+_i \otimes 1 + K_i \otimes X^+_i, \\
\De(X^-_i) = X^-_i \otimes K_i^{-1} + 1 \otimes X^-_i
\end{multline*}
and
\[
S(K_i) = K^{-1}_i, \,
S(X^+_i)= - K_i^{-1} X^+_i, \,
S(X^-_i)= - X^-_i K_i, \,
\epsilon (K_i)=1, \, \epsilon(X^\pm_i)=0.
\]
We refer to \cite[Ch. 4]{Ja} for details on the form of $\UU_q(\g)$ 
with this comultiplication.

Denote by $P$ and $P^+$ the sets of integral and dominant integral 
weights of $\g$. The sets of simple roots, simple coroots, and 
fundamental weights of $\g$ will be denoted by $\{\al_i\}_{i=1}^r,$ 
$\{\al_i\spcheck\}_{i=1}^r,$ and $\{\om_i\}_{i=1}^r,$ respectively.
Set $Q = \sum_{i=1}^r \Zset \al_i$ and 
$Q^+ = \sum_{i=1}^r \Nset \al_i$. Let $\lcor.,. \rcor$ be the symmetric
bilinear form on $\Span_\Qset \{ \al_1, \ldots, \al_r \}$
such that 
$\lcor \al_i, \al_j \rcor = d_i c_{ij}$. Recall the standard partial 
order on $P$:
\begin{equation}
\label{po}
\mbox{for} \; \;  \mu_1, \mu_2 \in P, \; 
\mu_1 < \mu_2, \; \; 
\mbox{if and only if} \; \; 
\mu_2 - \mu_1 \in Q^+ \backslash \{0\}.
\end{equation} 

Let $H$ be the group generated by $\{K_i^{\pm 1}\}_{i=1}^r$, 
which consists of all group like elements of $\UU_q(\g)$.
The $q$-weight spaces of an $H$-module $V$ are defined by
\[
V_\la = \{ v \in V \mid K_i v = q^{ \lcor \la, \al_i \rcor} v, \; 
\forall i = 1, \ldots, r \}, \; \la \in P.
\]
A $\UU_q(\g)$-module is called a type one module if it is the 
sum of its $q$-weight spaces. Each finite dimensional type one 
$\UU_q(\g)$-module is completely reducible \cite[Theorem 5.17]{Ja},
see the remark on p. 85 of \cite{Ja} for the validity of this for 
general base fields
$\KK$ and $q \in \KK^*$ not a root of unity. The category of 
(left) finite dimensional type one $\UU_q(\g)$-modules is closed under taking 
tensor products and duals (where the latter are defined as 
left modules using the antipode of $\UU_q(\g)$).  
The irreducible modules in this category  
are parametrized by $P^+$, see \cite[Theorem 5.10]{Ja}. 
We will denote by $V(\la)$ the irreducible finite 
dimensional type one $\UU_q(\g)$-module 
of highest weight $\la \in P^+$. 

Let $W$ and $\B_\g$ denote the Weyl and braid groups 
corresponding to $\g$. Let $s_1, \ldots, s_r$ be the simple reflections 
of $W$ corresponding to the roots $\al_1, \ldots, \al_r$, and
$T_1, \ldots, T_r$ be the related standard generators of $\B_\g$. 
Recall that the canonical projection $\B_\g \to W_\g$
has a section $W_\g \to \B_\g$, $w \mapsto T_w$ such that for each
reduced expression $w = s_{i_1} \ldots s_{i_l}$,
$T_w = T_{i_1} \ldots T_{i_l}$. The Bruhat order 
on $W$ will be denoted by $\leq$. For $w \in W$, 
we set $W^{\leq w} = \{ y \in W \mid y \leq w \}$.
We will use the $\B_\g$-action on $\UU_q(\g)$ given by 
\begin{align*}
& T_i(X_i^+) = - X_i^- K_i, \; T_i(X_i^-) = - K_i^{-1} X_i^+, \;
T_i(K_j) = K_j K_i^{-c_{ij}}, \\
& T_i(X_j^+) = \sum_{k=0}^{- c_{ij}} (- q_i)^{-k} (X_i^+)^{(-c_{ij}-k)}
X_j^+ (X_i^+)^{(k)}, \; j \neq i, \\
& T_i(X_j^-) = \sum_{k=0}^{- c_{ij}} (- q_i)^k (X_i^-)^{(k)}
X_j^- (X_i^-)^{(-c_{ij}-k)}, \; j \neq i,
\end{align*}
where $(X_i^\pm)^{(n)} = X_i^\pm/[n]_{q_i}!$,
see \cite[\S 8.14]{Ja} and \cite[\S 37.1]{L}. The braid group 
$\B_\g$ acts on all finite dimensional type one $\UU_q(\g)$-modules
$V$ by 
\[
T_i(v) = \sum_{l, m, n}
(-1)^m q_i^{m- ln} (X^+_i)^{(l)} (X^-_i)^{(m)} (X^+_i)^{(n)} v, 
\quad v \in V_\mu, \mu \in P,
\]
where the sum is over all $l,m, n \in \Nset$ such that 
$-l+m-n = \lcor \mu, \al_i\spcheck \rcor$, see
\cite[\S 8.6]{Ja} and 
\cite[\S 5.2]{L}. These actions satisfy
\begin{equation}
\label{comp}
T_w ( x . v ) = (T_w x) . (T_w v) \; \; 
\mbox{and} \; \; 
T_w(V(\la)_\mu) = V(\la)_{w \mu}
\end{equation} 
for all $w \in W$, $x \in \UU_q(\g)$, 
$v \in V(\la)$, $\la \in P^+$, $\mu \in P$,
see \cite[eq. 8.14(1)]{Ja}.
\subsection{}
\label{2.2}
Let $\UU_\pm$ be the subalgebras of $\UU_q(\g)$
generated by $\{X^\pm_i\}_{i=1}^r$. 
For a reduced decomposition
\begin{equation}
\label{wdecomp}
w = s_{i_1} \ldots s_{i_l}
\end{equation}
of an element $w \in W$, define the roots
\begin{equation}
\label{beta}
\beta_1 = \al_{i_1}, \beta_2 = s_{i_1} (\al_{i_2}), 
\ldots, \beta_l = s_{i_1} \ldots s_{i_{l-1}} (\al_{i_l})
\end{equation}
and Lusztig's root vectors
\begin{equation}
X^{\pm}_{\beta_1} = X^{\pm}_{i_1}, 
X^{\pm}_{\beta_2} = T_{s_{i_1}} (X^\pm_{i_2}), 
\ldots, X^\pm_{\beta_l} = T_{s_{i_1}} \ldots T_{s_{i_{l-1}}} (X^{\pm}_{i_l}),
\label{rootv}
\end{equation}
see \cite[\S 39.3]{L} and \cite[\S 8.24]{Ja}. 
In \cite[Proposition 2.2]{DKP}
De Concini, Kac and Procesi proved that the subalgebra
$\UU_\pm^w$ of $\UU_\pm$ generated by $X^{\pm}_{\beta_j}$, $j=1, \ldots, k$
does not depend on the choice of a reduced decomposition of $w$ 
and that it has the PBW type basis
\begin{equation}
\label{vect}
(X^\pm_{\beta_l})^{n_k} \ldots (X^\pm_{\beta_1})^{n_1}, \; \; 
n_1, \ldots, n_k \in \Nset.
\end{equation}
The fact that the space spanned by the monomials \eqref{vect} does not 
depend on the choice of a reduced decomposition of $w$ was independently 
proved by Lusztig \cite[Proposition 40.2.1]{L}. 
The elements $X^\pm_{\beta_l}$ satisfy the Levendorskii--Soibelman straightening 
rule, see \eqref{LS} below. A consequence of this is that all algebras 
$\UU^w_\pm$ are iterated skew polynomial rings and 
therefore are domains. 
There is a unique involutive algebra automorphism 
$\om$ of $\UU_q(\g)$ defined by 
\[
\om(X^\pm_i) = X^\mp_i, \;  
\om(K_i) = K_i^{-1}, \; \; 
i = 1, \ldots, r.
\]
It satisfies 
\begin{equation}
\label{om-prop}
\om( T_i(x)) =
(- q_i)^{ \lcor \al_i\spcheck, \ga \rcor }  T_i( \om (x)),
\; \; \forall x \in (\UU_q(\g))_\ga, \ga \in Q,
i=1,\ldots, r,
\end{equation}
see \cite[eq. 8.14(9)]{Ja}. This implies that for all 
$w \in W$, $x \in (\UU_q(\g))_\ga$, $\ga \in Q$
we have $\om(T_w(x)) = t' T_w (\om(x))$ for some $t' \in \KK^*$.
Therefore for all $w \in W$, $\om$ induces an algebra 
isomorphism between $\UU^w_+$ and $\UU^w_-$.

Denote the $\KK$-torus $\Tset^r = (\KK^*)^{\times r}$.
We have the group embeddings:
\begin{equation}
\label{embed}
H \hra P \hra \Tset^r,
\end{equation}
where the first one is given by $K_i \mt \al_i$, $i=1, \ldots, r$.  
The second one is given by $\om_i \mt (1, \ldots, 1, q, 1, \ldots, 1)$,
where $q$ is in position $i$. For $\ga \in Q$ denote the character
of $\Tset^r$:
\[
t \mt
t^\ga = \prod_{i=1}^r t_i^{\lcor \ga, \om_i \rcor}, \quad 
t =(t_1, \ldots, t_r) \in \Tset^r.
\]
The algebra $\UU_q(\g)$ is
$Q$-graded by
\begin{equation}
\label{Qgr}
\deg X^\pm_i = \pm \al_i, \, \, \deg K_i = 0.
\end{equation}
The homogeneous components of $\UU_q(\g)$ with respect to this grading
will be denoted by $(\UU_q(\g))_\ga$, $\ga \in Q$. 
There is a rational $\Tset^r$-action on $\UU_q(\g)$ by 
algebra automorphisms given by
\begin{equation}
\label{Tract}
t. x = t^\ga x, \quad \; \; \mbox{for} \; \; 
x \in (\UU_q(\g))_\ga, \ga \in Q.
\end{equation}
The conjugation action of $H$ on $\UU_q(\g)$ coincides
with the pull back of the action \eqref{Tract} under the 
embedding \eqref{embed}. The action \eqref{Tract} restricts 
to the following $P$-action on $\UU_q(\g)$:
\begin{equation}
\label{Pact} 
\mu \cdot x = q^{\lcor \mu, \ga \rcor} x, \quad \mbox{for}
\; \; x \in (\UU_q(\g))_\ga, \ga \in Q.
\end{equation}
The actions of $H$, $P$ and $\Tset^r$ on $\UU_q(\g)$ 
preserve the subalgebras $\UU^w_\pm$. 
\subsection{}
\label{2.3}
The quantum function algebra $R_q[G]$ is the Hopf  
subalgebra of the restricted dual of $\UU_q(\g)$ spanned by 
all matrix coefficients of the finite dimensional type one $\UU_q(\g)$-modules.
We think of $G$ as of the connected, simply connected algebraic group 
with Lie algebra $\g$, but $G$ is only used as a symbol 
since the base field $\KK$ is arbitrary 
(except that it cannot be finite, because it is assumed that 
$q \in \KK^*$ is not a root of unity). For $\la \in P^+$
the matrix coefficient of $\xi \in V(\la)^*$
and $v \in V(\la)$ will be denoted by 
$c_{\xi, v}^\la$; that is $c_{\xi, v}^\la(u) = \xi( u . v)$,
$\forall u \in \UU_q(\g)$.
Let $R^+$ be the subalgebra of $R_q[G]$ spanned by 
the matrix coefficients $c_{\xi , v}^\la$, where $\la \in P^+$, 
$\xi \in V(\la)^*$ and $v \in V(\la)_\la$. 
There are two canonical $H$-actions on $R_q[G]$ by algebra 
automorphisms:
\begin{equation}
\label{Hact}
K_i . c_{\xi,v}^\la = q^{ \lcor \nu ,\al_i \rcor } 
c_{\xi,v}^\la, \quad
K_i . c_{\xi,v}^\la = q^{ \lcor \mu,\al_i \rcor } 
c_{\xi,v}^\la, \quad 
\mbox{for} \; \;
\xi \in (V(\la)^*)_\nu, v \in V(\la)_\mu.
\end{equation}
One has a related $P \times P$-grading of $R_q[G]$
\begin{equation}
\label{PPgrad}
c_{\xi,v}^\la \in (R_q[G])_{\nu, \mu}, \quad 
\mbox{for} \; \; 
\xi \in (V(\la)^*)_\nu, v \in V(\la)_\mu.
\end{equation}

Fix highest weight vectors $v_\la \in V(\la)_\la$, $\la \in P^+$. For 
simplicity of the notation denote
\begin{equation}
\label{simpl}
c^\la_\xi = c^\la_{\xi, v_\la} \quad
\mbox{for} \; \; \la \in P^+, \xi \in V(\la)^*.
\end{equation}

All weight spaces $V(\la)_{w \la}= T_w (V(\la)_\la)$ are one dimensional.
For $\la \in P^+$ and $w \in W$ define 
$\xi_{w, \la} \in (V(\la)^*)_{- w\la}$ such that
$\lcor \xi_{w, \la}, T_w v_\la \rcor =1$. Let
\begin{equation}
\label{e}
e^\la_w = c^\la_{\xi_{w, \la}} = c^\la_{\xi_{w,\la}, v_\la}.
\end{equation}
Then 
\begin{equation}
\label{Prod}
e^{\la_1}_w e^{\la_2}_w = e^{\la_1 + \la_2}_w, \; \; 
\forall \la_1, \la_2 \in P^+, 
\end{equation}
see \cite[eq. (2.18)]{Y4} for our particular normalization 
and \cite[\S 9.1.10]{J} in general. Joseph proved
\cite[Lemma 9.1.10]{J} that
\[
E^+_w = \{ e_w^\la \mid \la \in P^+ \}
\]
is an Ore subset of $R^+$. The two actions \eqref{Hact} of $H$ 
on $R_q[G]$ and the $P\times P$-grading \eqref{PPgrad} 
descend to the localization
\[
R^w = R^+[(E_w^+)^{-1}].
\]
The invariant subalgebra of $R^w$ with respect to the $H$-action 
induced from the second action in \eqref{Hact} 
will be denoted by $R_0^w$. It was introduced by 
Joseph \cite[\S 10.4.8]{J} and called the 
quantum translated Bruhat cell. Its spectrum was studied by Gorelik 
in \cite{G}. 
Given $\mu \in P$, we decompose it as $\mu = \la_+ - \la_-$
for some $\la_\pm \in P^+$ and define 
$e_w^\mu = e_w^{\la_+} (e_w^{\la_-})^{-1}$. This does not depend on 
the choice of $\la_\pm$ because of \eqref{Prod}. 
We have that
\begin{equation}
\label{Rw0}
R^w_0 = \{ c^\la_\xi e_w^{-\la} \mid 
\la \in P^+, \xi \in V(\la)^* \}
\end{equation}
(in particular we do not need to take span in the right hand 
side), since 
\begin{multline}
\label{span}
\forall \la_1, \la_2 \in P^+, \; 
\xi \in V(\la_1)^*, \quad c^{\la_1}_\xi e^{-\la_1}_w =
c^{\la_1+\la_2}_{\xi'} e^{-\la_1- \la_2}_w,
\\
\mbox{where} \; \; 
\xi' = (\xi \otimes \xi_{w, \la_2})|_{\UU_q(\g)  (v_{\la_2} \otimes v_{\la_1} ) }
\in V(\la_1+\la_2)^*.
\end{multline}
Note that 
\begin{equation}
\label{grRw0}
R^w_0 = \bigoplus_{\ga \in Q} (R^w_0)_{\ga, 0}
\end{equation}
in terms of the induced grading from \eqref{PPgrad} and
\begin{equation}
\label{gr}
c^{\la}_\xi e^{-\la}_w \in (R^w_0)_{\nu + w(\la), 0}, \quad
\forall \xi \in (V(\la)^*)_\nu, \la \in P^+, \nu \in P. 
\end{equation} 

For $y \in W$ define the ideals
\begin{equation}
\label{id0}
Q(y)^\pm = \Span \{ c^\la_\xi \mid \la \in P^+, \, \xi \in V(\la)^*, \,
\xi \perp \UU_\pm T_y v_\la \}
\end{equation}
of $R^+$ and the ideals
\begin{equation}
\label{id2}
Q(y)^\pm_w = \{
c^\la_\xi e^{-\la}_w
\mid \la \in P^+, 
\xi \in V(\la)^*, \,
\xi \perp \UU_\pm T_y v_\la \}
\end{equation}
of $R^w_0$, see \cite{J,G} for details.
In the setting of \eqref{span}, one easily verifies that
$\xi \perp \UU_\pm T_y v_{\la_1}$ implies 
$\xi' \perp \UU_\pm T_y (v_{\la_2} \otimes v_{\la_2})$. Because of this,
there is no need to take span in the right hand side of \eqref{id2}.

\bth{Gorelik} (Gorelik) \cite[Lemmas 6.6 and 6.10, Corollary 7.1.2]{G}  
For all base fields $\KK$, $q \in \KK^*$ not a root of unity and 
a Weyl group element $w$ we have:

(a) The $H$-invariant prime ideals of the quantum translated 
Bruhat cell algebra $R_0^w$ (with respect to the first action 
\eqref{Hact}) which contain the ideal $Q(w)^+_w$ are the ideals
\[
Q(y)^-_w + Q(w)^+_w
\]
for $y \in W^{\leq w}$. All such ideals are completely prime.

(b) The poset of such ideals of $R_0^w$ ordered under inclusion,
is isomorphic to $W^{\leq w}$ equipped with the Bruhat order, i.e.
\[
Q(y_1)^-_w + Q(w)^+_w \subseteq Q(y_2)^-_w + Q(w)^+_w,
\]
if and only if $y_1 \leq y_2$. (In particular, all such 
ideals are distinct.)
\eth

Although we will not need this here, we note that Gorelik also 
described in \cite{G} all $H$-invariant prime ideals of $R^w_0$ (with respect 
to the first action \eqref{Hact}) in terms of the ideals $Q(y)^\pm_w$.
Gorelik stated the above results under the assumption that $\KK$ has 
characteristic $0$ and $q$ is transcendental over $\Qset$. However, 
her proofs work in the more general case when $q \in \KK^*$ is not a root of unity,
without any restrictions on $\KK$, see \cite[\S 3.2-3.4]{Y}.
\subsection{}
\label{2.4}
The quantum $R$-matrix associated to $w \in W$ is 
defined by
\begin{equation}
\RR^w = \prod_{j= k, \ldots, 1} \exp_{q_{i_j}}
\left( (q_{i_j}^{-1} - q_{i_j})
X^+_{\beta_j} \otimes X^-_{\beta_j} \right)
\label{Rw}
\end{equation}
in terms of Lusztig's root vectors \eqref{rootv}.
In \eqref{Rw} the noncommuting factors are multiplied in the 
order $j = k, \ldots, 1$, see e.g. \cite[eqs. 8.30(1) and 8.30(2)]{Ja}.
The $q$-exponential function is given by
\[
\exp_{q_i}(y) = \sum_{n=0}^\infty q_i^{-n(n-1)/2} 
\frac{y^n}{[n]_{q_i}!} \cdot
\]
The $R$-matrix $\RR^w$ belongs to a completion 
of $\UU^w_+ \otimes \UU^w_-$ and does not depend on the 
choice of a reduced decomposition of $w$, 
see \cite[\S 4.1.1]{L}.

There is a unique graded algebra antiautomorphism 
$\tau$ of $\UU_q(\g)$ defined by
\begin{equation}
\label{tau}
\tau(X_i^\pm) = X_i^\pm, 
\, 
\tau(K_i) = K_i^{-1}, \; \; 
i = 1, \ldots, r,
\end{equation}
cf. \cite[Lemma 4.6(b)]{Ja}. It satisfies
\begin{equation}
\label{tau-ident}
\tau (T_w x) = T_{w^{-1}}^{-1} ( \tau (x)), \; \; 
\forall w \in W, x \in \UU_q(\g), 
\end{equation}
see \cite[ eq. 8.18(6)]{Ja}. 

We will need the following result from \cite{Y,Y4}.

\bth{hom1} \cite[Theorem 2.6]{Y4} For all 
base fields $\KK$, $q \in \KK^*$ not a root of unity,
simple Lie algebras $\g$ and $w \in W$, the map
\[
\phi_w \colon R_0^w \to \UU^w_-, \; \; 
\phi_w(c^\la_\xi e_w^{-\la}) = 
(c^\la_{\xi, T_w v_\la} \otimes \id) (\tau \otimes \id) (\RR^w), \; 
\la \in P^+, \xi \in V(\la)^* 
\]
is a (well defined) surjective algebra antihomomorphism.
It is $H$-equivariant with respect to the first action \eqref{Hact}
of $H$ on $R_0^w$ and the conjugation action of $H$ on $\UU^w_-$.
The kernel of $\phi_w$ is $Q(w)^+_w$.
\eth
In the definition of $\phi_w$ the elements of $R_q[G]$ 
are viewed as functionals on $\UU_q(\g)$. The $H$-equivariance 
property is equivalent to saying that $\phi_w$ is graded, 
namely that $\phi_w( (R_0^w)_{\ga, 0}) = (\UU^w_-)_\ga$,
$\forall \ga \in Q$, cf. \eqref{Qgr}, \eqref{PPgrad}
and \eqref{grRw0}. A version of this theorem for $\UU_q(\g)$ 
equipped with the opposite comultiplication and for different 
braid group action and choice of Lusztig's root vectors  
was established in \cite[Theorem 3.7]{Y}.
To prove that such a map $\phi_w$ is well defined and is an 
algebra antihomomorphism, in \cite{Y,Y4} we first defined it
in terms of module algebras for Hopf algebras and then proved 
that it takes the above form.

We will use an interpretation of the algebras $\UU^w_-$ 
as quantized algebras of functions on Schubert cells 
using matrix coefficients of Demazure modules
from \cite{Y}. For $\la \in P^+$, $w \in W$
consider the Demazure modules 
$V_w(\la) = \UU_+ T_w v_\la = \UU^w_+ T_w v_\la$, 
see \cite[\S 4.4 and \S 6.3]{J} for details.
For $\eta \in V_w(\la)^*$ define 
\[
d^{w,\la}_{\eta} \in (\UU_+)^*, \; \; d^{w,\la}_{\eta}(x) = \lcor \eta, 
x T_w v_\la \rcor, \; x \in \UU_+.
\]
Set $U^w_+ = U_+ \cap w U_- w^{-1}$ where $U_\pm \subset G$ 
are the unipotent radicals of a pair of opposite Borel 
subgroups of the connected, simply connected algebraic group $G$
with Lie algebra $\g$. (We need those just as symbols
for reference purposes to the needed quantized coordinate 
rings, defined over an arbitrary base field.)
Denote 
by $R_q[U^w_+]$ the subset of $(\UU_+)^*$ consisting of 
\[
d^{w, \la}_\eta, \; \; \la \in P^+, \eta \in V_w(\la)^*.
\]
It is a $\KK$-vector space because of \eqref{span}.
In \cite[\S 3.8]{Y} we proved that 
\begin{equation}
\label{Rqmult}
d^{w, \la_1}_{\eta_1} d^{w, \la_2}_{\eta_2} 
= q^{\lcor \la_1, \la_1 + w^{-1}(\nu_1) \rcor}
d^{w, \la_1 + \la_2}_{\eta},
\end{equation}
where
\[
\eta : = \eta_1 \otimes \eta_2 |_{ \UU_+ (T_w v_{\la_1} \otimes T_w v_{\la_2})}
\; \; 
\mbox{and} \; \; 
\eta_i \in (V_w(\la_i)^*)_{\nu_i}
\]
defines an algebra structure on $R_q[U^w_+]$.
In particular, \eqref{Rqmult} is a well defined 
multiplication in $R_q[U^w_+]$. (The result in 
\cite[\S 3.8]{Y} concerned the Hopf algebra 
$\UU_q(\g)$ equipped with the opposite 
comultiplication and because of this there is a 
small difference in the power of $q$. The proofs are 
the same in both cases.) We have:

\bth{hom2} \cite[\S 3.8]{Y} For an arbitrary base field $\KK$, 
$q \in \KK^*$ not a root of unity, a simple Lie algebra $\g$, 
and a Weyl group element $w \in W$, we have: 

(a) The map
\begin{equation}
\label{vphi}
\varphi_w \colon R^w_0 \to R_q[U^w_+], \; \; 
\varphi_w ( c^\la_\xi e^{-\la}_w) = d^{w, \la}_{\xi |_{V_w(\la)} }, 
\; 
\la \in P^+, \xi \in V(\la)^*
\end{equation}
is a (well defined) surjective algebra homomorphism with kernel $Q(w)_w^+$.

(b) The algebras $R_q[U^w_+]$ and $\UU^w_-$ are antiisomorphic 
with an antiisomorphism given by 
\[
\psi_w \colon R_q[U^w_+] \to \UU^w_-, \quad
\psi_w ( d^{w, \la}_\eta ) =
(d^{w, \la}_\eta \otimes \id) 
(\tau \otimes \id)
(\RR^w), \; 
\la \in P^+, \eta \in V_w(\la)^*.
\]
\eth
 
We note that $\phi_w = \psi_w \varphi_w$.
In \cite[\S 3.8]{Y} we established the analog of \thref{hom2}
for the Hopf algebra $\UU_q(\g)$ equipped with the opposite 
comultiplication to the one considered here. \thref{hom2} 
follows from \thref{hom1} along the lines of the same 
argument as in \cite[\S 3.8]{Y}.

Since $Q(w)^+_w = \ker \varphi_w$ is a graded ideal of $R^w_0$ 
with respect to the $Q$-grading \eqref{grRw0} of $R^w_0$, 
one can push forward under $\varphi_w$ this grading 
to a $Q$-grading on $R_q[U^w_+]$. Comparing \eqref{gr} and \eqref{vphi}
gives that
\begin{equation}
\label{grRqU}
d^{w, \la}_\eta \in (R_q[U^w_+])_{\nu + w(\la)} 
\quad \mbox{for} \; \; 
\eta \in (V_w(\la)^*)_\nu, 
\end{equation}
where the weight spaces are computed with respect to the 
the action of $H$ on the dual of the Demazure module via 
the antipode of $\UU_q(\g)$. Since 
$\phi_w = \psi_w \varphi_w$ and $\phi_w$ 
is an antihomomorphism of $Q$-graded algebras, we 
have that $\psi_w$ is a graded antiisomorphism
with respect to the $Q$-gradings \eqref{grRqU} and 
\eqref{Qgr} of $R_q[U^w_+]$ and $\UU^w_-$, 
respectively.

It follows from \thref{hom2} (a) that the algebras 
$\UU^w_-$ are antiisomorphic to Joseph's algebras $S^+_w$, which are defined 
as the invariant subalgebras of the algebras $(R^+/Q(w)^+)[(E_w^+)^{-1}]$ 
with respect to the $H$-action induced from the second action 
in \eqref{Hact} (i.e. $S^+_w \cong R^w_0/Q(w)_w^+$).
We refer the reader to \cite[\S 10.3.1]{J} for 
details. These algebras played a key role in Joseph's 
work \cite{J1,J} on the spectrum of $R_q[G]$.

As it is customary in the area, here and below we denote
by the same symbols the images of elements of $R_q[G]$ and $R^+$
in their various quotients.    
\sectionnew{Polynormal generating sets of the $\Tset^r$-primes of $\UU_-^w$} 
\label{hsp}
\subsection{}
\label{3.1}
In this section we construct explicit $P$-polynormal 
generating sequences for all $\Tset^r$-prime ideals of 
the algebras $\UU^w_-$. In the case when the base field 
has characteristic $0$ and $q$ is transcendental over 
$\Qset$, the polynormal generating sets are very small,
see \thref{polynormal}. They correspond to certain subsets 
of the weight vectors of the duals of the Demazure modules of 
$\UU_q(\g)$ corresponding to fundamental weights. 
In the general case of an arbitrary base field $\KK$ 
and $q \in \KK^*$ not a root of unity, in \thref{pol2} 
we construct $P$-polynormal generating sets for all 
$\Tset^r$-invariant prime ideals of $\UU^w_-$. 
The generating sets correspond to (possibly) 
bigger sets derived from Demazure modules
for other highest weights. The proofs of Theorems 
\ref{tpolynormal} and \ref{tpol2} are analogous.

For the convenience of the reader we recall several definitions 
regarding polynormality. Assume that $I$ is an ideal 
of a ring $R$. A sequence $u_1, \ldots, u_n \in R$ 
is called a polynormal generating sequence if the set
$\{u_1, \ldots, u_n\}$
generates $I$ and for all $i = 1, \ldots, n$, the element $u_i$
is normal in $R$ modulo the ideal generated by $u_1, \ldots, u_{i-1}$. 
In particular $u_1$ should be a normal element of $R$. 
If $u_1, \ldots, u_n$ is a polynormal generating sequence 
of $I$, then $I$ is generated both as a left and right ideal of $R$
by the set $\{u_1, \ldots, u_n \}$. 

If a group $\Ga$ acts on the ring $R$ by algebra automorphisms, 
we say that an element $u \in R$ is $\Ga$-normal if it is a
$\Ga$-eigenvector and if there exists $g \in \Ga$ such that
$u r = (g . r) u$ for all $r \in R$. We note that sometimes 
$\Ga$-normality is defined requiring only the second condition, 
see \cite{G1}. In all cases we will be able to construct 
elements satisfying both conditions. We also note that 
requiring only the second condition will not be sufficient 
to extend this definition to $\Ga$-polynormality, as we do next.
We say that an element $u \in R$ is $\Ga$-normal modulo a 
$\Ga$-stable ideal $I$, if its image in $R/I$ is $\Ga$-normal.

We say that a sequence $u_1, \ldots, u_n \in R$ is a 
$\Ga$-polynormal generating sequence of a $\Ga$-stable ideal $I$ if 
$\{u_1, \ldots, u_n\}$ generates $I$ and for all $i = 1, \ldots, n$,
the element $u_i$ is a $\Ga$-normal element of $R$ modulo
the ideal generated by $u_1, \ldots, u_{i-1}$. We note that 
the conditions posed on the elements $u_1, \ldots, u_{i-1}$
imply that the ideal of $R$ generated by them is $\Ga$-stable.

For $y, w \in W$, $y \leq w$, define
\begin{multline}
\label{Iw}
I_w(y) = \phi_w (Q(y)_w^-+Q(w)_w^+) = \phi_w (Q(y)_w^-) \\
= \{ (d^{w,\la}_\eta \otimes \id)(\RR^w) \mid 
\la \in P^+, \eta \in (V_w(\la) \cap \UU_- T_y v_\la)^\perp \}
\subset \UU^w_-.
\end{multline}
\thref{Gorelik} of Gorelik \cite{G} and \thref{hom1} imply the first 
two parts of the following theorem. The third part of the 
theorem is \cite[Theorem 1.1]{Y}.

\bth{Yold} For an arbitrary base field $\KK$, $q \in \KK^*$ not 
a root of unity, a simple Lie algebra $\g$, 
and a Weyl group element $w \in W$, we have: 
 
(a) If $y \in W^{\leq w}$, then $I_w(y)$ is a $\Tset^r$-invariant completely 
prime ideal of $\UU^w_-$ with respect to the action \eqref{Tract}. 
All $\Tset^r$-invariant prime ideals of $\UU^w_-$ are of this form.

(b) The correspondence $y \in W^{\leq w} \mapsto I_w(y)$ is an isomorphism
from the poset $W^{\leq w}$ equipped with the Bruhat order
to the poset of $\Tset^r$-invariant 
prime ideals of $\UU^w_-$ ordered under inclusion; that is
$I_w(y) \subseteq I_w(y')$ for $y, y' \in W^{\leq w}$
if and only if $y \leq y'$.

(c) Assume that $\KK$ has characteristic $0$ and $q$ is transcendental 
over $\Qset$. Then $I_w(y)$ is generated as a right ideal by
\[
\psi_w(d^{w, \om_i}_\eta)= (d^{w, \om_i}_\eta \otimes \id) (\RR^w) \quad 
\mbox{for} \quad
\eta \in (V_w(\om_i) \cap \UU_- T_y v_{\om_i})^\perp, i= 1, \ldots, r,
\]
where $\om_1, \ldots, \om_r$ are the fundamental weights 
of $\g$.
\eth

For part (c) one needs the stronger assumptions on $q$ 
and $\KK$ (that $\KK$ has characteristic 0 and $q$ is 
transcendental over $\Qset$), because it
relies on Joseph's result \cite[Th\'eor\`eme 3]{J2} which uses 
a specialization argument.

Note that the sets of invariant subspaces of $\UU^w_-$ with 
respect to the conjugation action of $H$, the $P$-action \eqref{Pact}, 
and the $\Tset^r$-action \eqref{Tract} coincide. Thus 
\begin{equation}
\label{HPT}
H - \Spec \UU^w_- = P - \Spec \UU^w_- = \Tset^r - \Spec \UU^w_-.
\end{equation}
We use $\Tset^r$-invariance in \thref{Yold} to align our treatment 
to the Goodearl--Letzter framework \cite{GL}. In \S\ref{2.3}-\ref{2.4} 
we used $H$-invariance instead, because it was more convenient 
to state the results within the framework of the adjoint 
action of the Hopf algebra $\UU_q(\g)$ on itself,
although an appropriate torus invariance 
could have been used as well.
\subsection{}
\label{3.2}
Denote by $w_0$ the longest element of $W$ and 
set $\RR = \RR^{w_0}$. For $\ga \in Q^+$, 
$\ga \neq 0$ denote 
$m(\ga) = \dim (\UU_+)_\ga= \dim (\UU_-)_{-\ga}$,
and fix a pair of dual bases 
$\{u_{\ga, k} \}_{k=1}^{m(\ga)}$ and
$\{u_{-\ga, k} \}_{k=1}^{m(\ga)}$
of $(\UU_+)_\ga$ and $(\UU_-)_{-\ga}$ 
with respect to the Rosso--Tanisaki form,
see \cite[Ch. 6]{Ja}. Then
\begin{equation}
\label{Rm}
\RR = 1 \otimes 1 + \sum_{\ga \in Q^+, \ga \neq 0} 
\sum_{k=1}^{m(\ga)} u_{\ga, k} \otimes u_{-\ga, k}. 
\end{equation}
Recall the standard $R$-matrix commutation relations:

\ble{comm}
For all 
$\la_i \in P^+$, $\nu_i \in P$, 
$\xi_i \in V(\la_i)^*_{\nu_i}$, $i=1,2$:   
\begin{multline*}
c_{\xi_1}^{\la_1} c_{\xi_2}^{\la_2} =
q^{ \lcor \la_1, \la_2 \rcor - \lcor \nu_1, \nu_2 \rcor} 
c_{\xi_2}^{\la_2} c_{\xi_1}^{\la_1} 
\\
+ 
\sum_{\ga \in Q^+, \ga \neq 0}
\sum_{k=1}^{m(\ga)}
q^{ \lcor \la_1, \la_2 \rcor - \lcor \nu_1 - \ga , \nu_2 + \ga \rcor}
c_{S^{-1}(u_{\ga, k})\xi_2}^{\la_2} 
c_{S^{-1}(u_{-\ga, k}) \xi_1}^{\la_1}.
\end{multline*}
\ele

For details we refer to \cite[Theorem I.8.15]{BG}. 
\leref{comm} implies that for all $\la_1, \la_2 \in P^+$,
$\nu_2 \in P$, $\xi_2 \in V(\la_2)^*_{\nu_2}$
\begin{equation}
\label{cw-comm}
c_{\xi_1}^{\la_1} e_w^{\la_2} = 
q^{\lcor \la_1 , \la_2 \rcor + \lcor \nu_1, w(\la_2) \rcor} 
e_w^{\la_2} c_{\xi_1}^{\la_1} \mod Q(w)^+,
\end{equation}
recall \eqref{id0}.
Combining \leref{comm} and \eqref{cw-comm} leads to the following 
result.

\ble{comm-main}
For all 
$\la_i \in P^+$, $\nu_i \in P$, 
$\xi_i \in V(\la_i)^*_{\nu_i}$, $i=1,2$
\begin{multline*}
( c_{\xi_1}^{\la_1} e^{-\la_1}_w ) 
( c_{\xi_2}^{\la_2} e^{-\la_2}_w ) -
q^{ \lcor w(\la_1) - \nu_1 , w(\la_2) + \nu_2 \rcor} 
( c_{\xi_2}^{\la_2} e^{-\la_2}_w ) 
( c_{\xi_1}^{\la_1} e^{-\la_1}_w ) -
\\
\sum_{\ga \in Q^+, \ga \neq 0}
\sum_{k=1}^{m(\ga)}
q^{ 
\lcor w(\la_1) - \nu_1 + \ga , w(\la_2) + \nu_2 + \ga \rcor - 
\lcor w (\la_1), \ga \rcor}
( c_{S^{-1}(u_{\ga, k})\xi_2}^{\la_2} e^{-\la_2}_w ) 
( c_{S^{-1}(u_{-\ga, k}) \xi_1}^{\la_1} e^{-\la_1}_w )
\end{multline*}
belongs to $Q(w)^+_w$.
\ele
\subsection{}
\label{3.3}
Fix $y \in W^{\leq w}$. For each $i=1, \ldots, r$ choose a
basis $\Om_i$ of the orthogonal complement 
$(V_w(\om_i) \cap \UU_- T_y v_{\om_i})^\perp$ inside $V_w(\om_i)^*$,
consisting of weight vectors (with respect to the $H$-action).
Let $\Om_w(y) = \Om_1 \sqcup \ldots \sqcup \Om_r$. For 
$\eta \in (V_w(\om_i))^*_\nu$ denote $\iota(\eta) = i$ 
and $\Wt(\eta) = \nu$. In particular, we 
have maps $\iota \colon \Om_w(y) \to \{1, \ldots, r\}$
and $\Wt \colon \Om_w(y) \to P$. 

\bth{polynormal} Assume that $\KK$ has characteristic $0$ and $q$ is 
transcendental over $\Qset$.
Choose any linear ordering on $\Om_w(y)$ 
with the property that if $\eta_1, \eta_2 \in \Om_w(y)$, 
$\iota(\eta_1) = \iota(\eta_2)$ and 
$\Wt(\eta_1) < \Wt(\eta_2)$, then 
$\eta_1 < \eta_2$ (recall \eqref{po}).
Denote $\Om_w(y)= \{ \eta_1 < \eta_2 < \ldots < \eta_{|\Om_w(y)|} \}$.
Then 
\begin{equation}
\label{gen}
\psi_w \Big( d_{\eta_j}^{w, \om_{\iota(\eta_j)}} \Big)=
\Big( d^{w, \om_{\iota(\eta_j)}  }_{\eta_j} \otimes \id \Big) 
(\RR^w), \quad j = 1, \ldots, |\Om_w(y)| 
\end{equation}
is a $P$-polynormal generating sequence for the $\Tset^r$-prime 
ideal $I_w(y)$ of $\UU^w_-$, with respect to the action \eqref{Pact}.
\eth
In \cite{Cal} Caldero obtained a nonconstructive proof of the polynormality 
in the case $w= w_0$  (the longest element of the Weyl group $W$) 
and $\KK = \Cset(q)$.
\medskip
\\
{\em{Proof of \thref{polynormal}}}.
\thref{Yold} (c) implies that the set \eqref{gen} 
generates the ideal $I_w(y)$. Its elements are homogeneous with 
respect to the grading \eqref{Qgr} and are thus $P$-eigenvectors
with respect to the action \eqref{Pact}. More precisely
\[
\psi_w \Big( d_{\eta_j}^{w, \om_{\iota(\eta_j)}} \Big)
\in (\UU^w_-)_{\Wt(\eta_j) + w(\om_{\iota(\eta_j)})} , \quad 
j = 1, \ldots, |\Om_w(y)|,
\]
cf. \eqref{grRqU}. Denote $i_j = \iota(\eta_j)$
and $\nu_j = \Wt(\eta_j)$. Since
\[
\eta_j \in (V_w(\om_{i_j}))^*_{\nu_j} \; \; 
\mbox{and} \; \;  
\eta_j \perp (V_w(\om_{i_j}) \cap \UU_- T_y v_{\om_{i_j}}),
\]
there exist preimages $\xi_j \in V(\om_{i_j})^*$ with the following 
properties
\[
\xi_j \in (V(\om_{i_j}))^*_{\nu_j}, \;
\xi_j \perp \UU_- T_y v_{\om_{i_j}} \; \;
\mbox{and} \; \; 
\xi_j|_{V_w(\om_{i_j})} = \eta_j. 
\]
We fix a family of such preimages.
Applying \thref{hom2} and \leref{comm-main}, we obtain that for all 
$\xi \in (V(\la))^*_\nu$, $\la \in P^+$, $\nu \in P$:
\begin{align*}
&\phi_w(c_{\xi}^{\la} e^{-\la}_w )
\psi_w \Big( d_{\eta_j}^{w, \om_{i_j}} \Big) 
-
\psi_w \Big( d_{\eta_j}^{w, \om_{i_j}} \Big)
\Big( ( w(\om_{i_j}) - \nu_j) \cdot 
\phi_w( c_{\xi}^{\la} e^{-\la}_w ) \Big)
\\ 
=&
\sum_{\ga \in Q^+, \ga \neq 0}
\sum_{k=1}^{m(\ga)}
q^{ \lcor w(\om_{i_j}) - \nu_j + \ga , w(\la) + \nu + \ga \rcor
- \lcor w(\om_{i_j}), \ga \rcor }
\phi_w \Big(  
c_{ S^{-1}(u_{-\ga, k}) \xi_j}^{ \om_{i_j} } e^{-\om_{i_j}}_w \Big)
\\
& \hspace{8.5cm} \times
\phi_w(c_{S^{-1}(u_{\ga, k})\xi}^{\la} e^{-\la}_w)
\end{align*}
in terms of the $P$-action \eqref{Pact} and the dual bases 
$\{ u_{\pm \ga, k } \}_{k=1}^{m(\ga)}$ from \S \ref{3.2}.
Recall from \thref{hom1} that 
$\phi_w \colon R^w_0 \to \UU^w_-$ is a graded antihomomorphism. Thus 
$\phi_w( c^\la_\xi e^{-\la}_w ) \in 
(\UU^w_-)_{\nu + w(\la)}$ and 
\[
(w(\om_{i_j})-\nu_j) \cdot \phi_w ( c^\la_\xi e^{-\la}_w )
= q^{ \lcor w(\om_{i_j})-\nu_j, w(\la) + \nu \rcor } 
\phi_w ( c^\la_\xi e^{-\la}_w ).
\]
Since $\xi_j \in (\UU_- T_y v_{\om_{i_j}})^\perp$ 
and $(\UU_- T_y v_{\om_{i_j}})^\perp$
is a $\UU_-$-submodule of $V(\om_{i_j})^*$, we have that 
\[
S^{-1}(u_{-\ga, k}) \xi_j \in ( \UU_- T_y v_{\om_{i_j}} )^\perp,
\]
for all $\ga \in Q^+ \backslash \{0 \}$ and $k = 1, \ldots, m(w)$.  
The order relation on $\Om_w(y)$ implies that for all $\ga' \in Q^+$, 
$\ga' \neq 0$
\[
( (\UU_- T_y v_{\om_{i_j}})^\perp )_{\nu_j - \ga'} \subseteq 
\Span \{ \xi_n \mid n=1, \ldots, j-1, \; 
\iota(n) = i_j \} +
V_w( {\om_{i_j}} )^\perp.
\]
Therefore
\[
S^{-1}(u_{-\ga, k}) \xi_j
\in \Span \{ \xi_n \mid n=1, \ldots, j-1, \; 
\iota(n) = i_j \} +
V_w( {\om_{i_j}} )^\perp
\]
and
\begin{multline*}
\phi_w \Big( 
c_{ S^{-1}(u_{-\ga, k}) \xi_j }^{ \om_{i_j} } e^{-\om_{i_j}}_w
\Big)
\in 
\Big\lcor 
\phi_w \Big( c^{\om_{i_1}}_{\xi_1} e_w^{-\om_{i_1}} \Big),
\ldots, 
\phi_w \Big( c^{\om_{i_{j-1}}}_{\xi_{j-1} } e_w^{-\om_{i_{j-1}}} \Big)
\Big\rcor
\\
=
\Big\lcor 
\psi_w \Big( d_{\eta_1}^{w, \om_{i_1}} \Big),
\ldots, 
\psi_w \Big( d_{\eta_{j-1}}^{w, \om_{i_{j-1}}} \Big)
\Big\rcor
\end{multline*}
for all $\ga \in Q^+$, $\ga \neq 0$ and $k = 1, \ldots, m(\ga)$.
Therefore 
\[
\psi_w \Big( d_{\eta_j}^{w, \om_{i_j}} \Big) u
- ( (\nu_j - w(\om_{i_j})) \cdot u )  
\psi_w \Big( d_{\eta_j}^{w, \om_{i_j}} \Big)
\in \Big\lcor 
\psi_w \Big( d_{\eta_1}^{w, \om_{i_1}} \Big),
\ldots, 
\psi_w \Big( d_{\eta_{j-1}}^{w, \om_{i_{j-1}}} \Big)
\Big\rcor
\]
for all $u \in \UU^w_-$ which proves the statement of the 
theorem.
\qed
\bre{rem} Because of the embedding \eqref{embed}, \thref{polynormal} (and \thref{pol2} below)
are stronger results than constructing $\Tset^r$-polynormal generating sequences 
for the $\Tset^r$-prime ideals of $\UU^w_-$. Recalling the definition of 
$P$-polynormality and $\Tset^r$-polynormality \S \ref{3.1}, one should note 
that the $P$-eigenvectors and $\Tset^r$-eigenvectors in $\UU^w_-$ are the 
same, as they are simply the homogeneous elements of $\UU^w_-$ with 
respect to the grading \eqref{Qgr}. 
\ere
\subsection{}
\label{3.4} 
Next we prove that all ideals $I_w(y)$ of $\UU^w_-$ are 
$P$-polynormal under the weaker assumption that $q \in \KK^*$ is not a root
of unity without any restrictions on the characteristic of the field 
$\KK$.

Recall \eqref{Iw}.
Since $\UU^w_-$ is noetherian, for each $y \in W^{\leq w}$ 
there exists a finite set $\Sig_w(y) \subset P^+$ such that
\[
\{ \psi_w( d^{w, \la}_\eta) \mid 
\la \in \Sig_w(y), 
\eta \in (V_w(\la) \cap \UU_- T_y v_\la)^\perp \subset V_w(\la)^* \}
\]
generates $I_w(y)$. For $\la \in \Sig_w(y)$ let 
$\Ga_\la$ be a basis of the orthogonal complement 
$(V_w(\la) \cap \UU_- T_y v_\la)^\perp$ inside $V_w(\la)^*$,
which consists of weight vectors (with respect to the 
action of $H$).
Denote $\Ga_w(y) = \sqcup \{ \Ga_\la \mid \la \in \Sig_w(y) \}$.
Define the maps 
\[
\Hw \colon \Ga_w(y) \to \Sig_w(y) \; \; 
\mbox{and} \; \; 
\Wt \colon \Ga_w(y) \to P
\]
by
\[
\Hw(\eta) = \la, \; \; \Wt(\eta) = \nu, \; \; 
\mbox{if} \; \; \eta \in (V_w(\la))^*_\nu,
\] 
where $\Hw(.)$ stands for highest weight and $\Wt(.)$ 
stands for weight.

\bth{pol2} Let $\KK$ be an arbitrary base field, $q \in \KK^*$
not a root of unity, and $y \in W$, $y \leq w$. 
Choose a set $\Sig_w(y)$ as above. Consider any linear 
ordering on $\Ga_w(y)$ 
with the property that if $\eta_1, \eta_2 \in \Ga_w(y)$, 
$\Hw(\eta_1) = \Hw(\eta_2)$ and 
$\Wt(\eta_1) < \Wt(\eta_2)$, then 
$\eta_1 < \eta_2$. Let
$\Ga_w(y)= \{ \eta_1 < \eta_2 < \ldots < \eta_{|\Ga_w(y)|} \}$.
Then 
\begin{equation}
\label{gen2}
\psi_w \Big( d_{\eta_j}^{w, \Hw(\eta_j)} \Big)=
\Big( d^{w, \Hw(\eta_j)  }_{\eta_j} \otimes \id \Big) 
(\RR^w), \quad j = 1, \ldots, |\Ga_w(y)| 
\end{equation}
is a $P$-polynormal generating sequence for the $\Tset^r$-prime 
ideal $I_w(y)$ of $\UU^w_-$. Moreover,
\begin{multline}
\label{normal}
\psi_w \Big( d_{\eta_j}^{w, \Hw(\eta_j)} \Big) u
=[ (\Wt(\eta_j) - w(\Hw(\eta_j) ) ) \cdot u ]
\psi_w \Big( d_{\eta_j}^{w, \Hw(\eta_j) } \Big)
\\
\mod
\Big\lcor 
\psi_w \Big( d_{\eta_1}^{w, \Hw(\eta_1) } \Big),
\ldots, 
\psi_w \Big( d_{\eta_{j-1}}^{w, \Hw(\eta_j) } \Big)
\Big\rcor
\end{multline}
for all $j =1, \ldots, |\Ga_w(y)|$, $u \in \UU^w_-$.
\eth
\begin{proof} We argue analogously to the proof of \thref{polynormal}.
The choice of the set $\Ga_w(y)$ implies that the elements \eqref{gen2}
generate the ideal $I_w(y)$. They are homogeneous with 
respect to the grading \eqref{Qgr} 
\[
\psi_w \Big(
d_{\eta_j}^{w, \Hw(\eta_j) }
\Big)
\in (\UU^w_-)_{\Wt(\eta_j) + w(\Hw(\eta_j))} , \quad 
j = 1, \ldots, |\Gamma_w(y)|,
\]
see \eqref{grRqU}, and are thus $P$-eigenvectors
with respect to the action \eqref{Pact}.

From the definition of the elements $\eta_j$ it follows that
there exist preimages $\xi_j \in V( \Hw(\eta_j) )^*$ such that
\[
\xi_j \in (V( \Hw(\eta_j) ))^*_{ \Wt(\eta_j) }, \;
\xi_j \perp \UU_- T_y v_{ \Hw(\eta_j) } \; \;
\mbox{and} \; \; 
\xi_j|_{V_w( \Hw(\eta_j) )} = \eta_j. 
\]
We fix a family of such preimages. \thref{hom2} and \leref{comm-main} 
imply that for all
for all $\xi \in (V(\la))^*_\nu$, $\la \in P^+$, $\nu \in P$:
\begin{align*}
&\phi_w(c_{\xi}^{\la} e^{-\la}_w ) 
\psi_w \Big( d_{\eta_j}^{w, \Hw(\eta_j) } \Big)
-
\psi_w \Big( d_{\eta_j}^{w, \Hw(\eta_j) } \Big)
\Big( (w(\Hw(\eta_j)) - \Wt(\eta_j) ) \cdot 
\phi_w( c_{\xi}^{\la} e^{-\la}_w ) \Big)
\\ 
=&
\sum_{\ga \in Q^+, \ga \neq 0}
\sum_{k=1}^{m(\ga)}
q^{ \lcor w(\Hw(\eta_j)) - \Wt(\eta_j) + \ga, w(\la) + \nu + \ga\rcor 
- \lcor w(\Hw(\eta_j)) , \ga \rcor }
\phi_w \Big( 
c_{ S^{-1}(u_{-\ga, k}) \xi_j}^{ \Hw(\eta_j) } 
e^{- \Hw(\eta_j) }_w \Big)
\\
&\hspace{8.5cm} \times
\phi_w(  c_{S^{-1}(u_{\ga, k})\xi}^{\la} e^{-\la}_w ),
\end{align*}
where $\{ u_{\pm \ga, k} \}_{k=1}^{m(\ga)}$ are 
dual bases of $(\UU_\pm)_{ \pm \ga}$ as in \S \ref{3.2}.
As in the proof of \thref{polynormal}, 
the properties of the linear ordering of $\Ga_w(y)$ imply that
\[
S^{-1}(u_{-\ga, k}) \xi_j
\in \Span \{ \xi_n \mid n=1, \ldots, j-1, \; 
\Hw(\eta_n) = \Hw(\eta_j) \} +
V_w( \Hw(\eta_j) )^\perp
\]
and 
\[
\phi_w \Big( 
c_{ S^{-1}(u_{-\ga, k}) \xi_j }^{ \Hw(\eta_j) } 
e^{- \Hw(\eta_j) }_w
\Big)
\in 
\Big\lcor 
\psi_w \Big( d_{\eta_1}^{w, \Hw(\eta_j)} \Big),
\ldots, 
\psi_w \Big( d_{\eta_{j-1}}^{w, \Hw(\eta_j) } \Big)
\Big\rcor
\]
for all $\ga \in Q^+$, $\ga \neq 0$ and $k = 1, \ldots, m(\ga)$.
This implies \eqref{normal}.
\end{proof} 
\sectionnew{The Goodearl--Lenagan conjecture} 
\label{GL}
\subsection{}
\label{4.1}
In this section
we specialize the results from the previous one
to obtain a constructive proof of the Goodearl--Lenagan 
conjecture \cite{GLen} that all torus invariant prime ideals
of the algebras of quantum matrices have polynormal 
generating sets consisting of quantum minors. Since all fundamental 
representations of ${\mathfrak{sl}}_{r+1}$ are minuscule, 
all ideal generators in \thref{polynormal} become quantum minors.

Fix two positive integers $m$ and $n$.
Recall that the algebra of quantum matrices $R_q[M_{m,n}]$ is the 
$\KK$-algebra with generators $x_{ij}$, $i= 1, \ldots,m$,
$j = 1, \ldots, n$ and relations
\begin{align*}
x_{ij} x_{lj} &= q x_{lj} x_{ij}, \quad \mbox{for} \; i < l, \\
x_{ij} x_{ik} &= q x_{ik} x_{ij}, \quad \mbox{for} \; j < k, \\ 
x_{ij} x_{lk} &= x_{lk} x_{ij}, \quad 
\mbox{for} \; i < l, j > k,\\
x_{ij} x_{lk} - x_{lk} x_{ij} &= (q-q^{-1}) x_{ik} x_{lj}, \quad 
\mbox{for} \; i < l, j<k.
\end{align*}
The torus $\Tset^{m+n} = (\KK^*)^{ \times (m+n)}$ acts on $R_q[\Mmn]$ 
by algebra automorphisms by 
\begin{equation}
\label{torus2}
(t_1, \ldots, t_{m+n}) \cdot x_{ij} 
= t_i t_{m+j}^{-1} x_{ij}, 
\quad (t_1, \ldots, t_{m+n}) \in \Tset^{m+n}.
\end{equation}
For two integers $i \leq l$ denote $[i,l] = \{i, i+1, \ldots, l \}$. 
Given two subsets $J = \{ j_1 < \ldots < j_k \} \subset [1,m]$
and $J' = \{ j'_1 < \ldots < j'_k \} \subset [1,n]$
define the quantum minor $\De^q_{J,J'} \in R_q[\Mmn]$
by 
\[
\De^q_{J,J'} = \sum_{w \in S_k} (-q)^{l(w)}
x_{j_1 j'_{w(1)}} \ldots x_{j_k j'_{w(k)}} 
= \sum_{w \in S_k} (-q)^{-l(w)}
x_{j_{w(k)} j'_k} \ldots x_{j_{w(1)} j'_1}. 
\]

Throughout this subsection we set $\g = \sl_{m+n}$, $W = S_{m+n}$, and 
$w= c^m$, where $c$ is the Coxeter element $(1 2 \ldots m+n)$.
M\'eriaux and Cauchon constructed an isomorphism between 
$\UU^{c^m}_+$ and $R_q[\Mmn]$ in \cite[Proposition 2.1.1]{MC}.
Recall from \S \ref{2.2} that the automorphism $\om$ of 
$\UU_q(\sl_{m+n})$ restricts to an isomorphism 
$\om \colon \UU^{c^m}_+ \to \UU^{c^m}_-$. Furthermore 
we have the antiisomorphism 
$\psi_{c^m} \colon R_q[U_+^{c^m}] \to \UU^{c^m}_-$
from \thref{hom2} (b).
The composition of the 
above maps provides an antiisomorphism 
$R_q[U_+^{c^m}] \to R_q[\Mmn]$, which was used 
in \cite{Y} to study $\Tset^{m+n} -\Spec R_q[\Mmn]$. We 
briefly go over the modifications of this construction,
which are needed because of the different comultiplication 
of $\UU_q(\g)$, braid group action and Lusztig's root vectors 
of $\UU_q(\g)$ used in this paper. 

Consider the reduced expression
\[
c^m = (s_m \ldots s_1) (s_{m+1} \ldots s_2) \ldots 
(s_{m+n-1} \ldots s_n).
\]
Denote the Lusztig root vectors of 
$\UU^{c^m}_+$ and $\UU^{c^m}_-$
constructed in \eqref{rootv} by 
\[
X_{1,m+1}, \ldots, X_{m,m+1}; 
X_{1,m+2}, \ldots, X_{m,m+2}; \ldots; 
X_{1,m+n}, \ldots, X_{m,m+n}.
\]
and
\[
X_{m+1,1}, \ldots, X_{m+1,m}; X_{m+2,1}, \ldots, X_{m+2,m}; \ldots; 
X_{m+n,1}, \ldots, X_{m+n,m},
\]
respectively. For $i \in [1,m]$ and $j \in [1,n]$ set 
\begin{equation}
\label{zetamn}
\zeta_{m,n}(X_{m+j, i}) = (-q)^{i+j-2} x_{ij}.
\end{equation}

\ble{Rqisom} \cite[Lemma 5.4]{Y} The map \eqref{zetamn}
extends (uniquely) to an algebra isomorphism 
$\zeta_{m,n} \colon \UU^{c^m}_- \to R_q[\Mmn]$.
\ele
\begin{proof} By \cite[Proposition 2.1.1]{MC}, the
map $X_{i,m+j} \mt x_{ij}$, $i \in [1,m]$, $j \in [1,n]$ 
defines an isomorphism $\UU^{c^m}_+ \cong R_q[\Mmn]$. 
By iterating \eqref{om-prop}, one obtains that the isomorphism 
$\om \colon \UU^{c^m}_+ \to \UU^{c^m}_-$ satisfies 
\[
\om(X_{i, m+j}) = (-q)^{2-i-j} X_{m+j,i}, \; \; 
\forall i \in [1,m], j \in [1,n].  
\]    
The map $\zeta_{m,n}$ is equal to an appropriate composition of 
these isomorphisms.
\end{proof}

For all $i \in [1,m]$, $j \in [1,n]$,
\begin{equation}
\label{Xcom}
X_{j+m,i} \in (\UU_q(\sl_{m+n}))_{- \al_{m-i+1} - \ldots - \al_{m+j-1}}.
\end{equation}
Because of this $\zeta_{m,n}$ defines a bijection 
between the set of $\Tset^{m+n-1}$-eigenvectors in 
$\UU^{c^m}_-$ with respect to the action \eqref{Tract} and the set 
of $\Tset^{m+n}$-eigenvectors in $R_q[\Mmn]$
with respect to \eqref{torus2}. Denote by $P$ 
the weight lattice of $\sl_{m+n}$. It follows from 
\eqref{Xcom} that the $P$-action \eqref{Pact}
on $\UU^{c^m}_-$ transfers via $\zeta_{m,n}$ to the following 
action of $P$ on $R_q[M_{m,n}]$:
\begin{equation}
\label{newPact}
\om_k \cdot x_{ij} = q^{ - \de_{k, [m-i+1, m+j-1] } } x_{ij}
\end{equation}
for $k \in [1,m+n-1]$, $i \in [1,m]$, $j \in [1,n]$,
where for $a,b,c \in \Zset$, $\delta_{a, [b,c]} = 1$ 
if $a \in [b,c]$ and $\delta_{a, [b,c]} = 0$ otherwise.

For $y \in S_{m+n}^{\leq c^m}$ denote
\begin{equation}
\label{ideal}
I_{m,n}(y) = \zeta_{m,n}(I_{c^m}(y)).
\end{equation}
\thref{Yold} implies:

\bco{RqMmn} Let $\KK$ be an arbitrary base field 
and $q \in \KK^*$ not a root of unity.
For all $y \in S_{m+n}^{\leq c^m}$, $I_{m,n}(y)$
is a $\Tset^{m+n}$-invariant prime ideal of $R_q[\Mmn]$ and all
$\Tset^{m+n}$-primes of $R_q[\Mmn]$ are of this form. The 
map $S_{m+n}^{\leq c^m} \to \Tset^{m+n} -\Spec R_q[\Mmn]$ 
given by $y \mapsto I_{m,n}(y)$ is an isomorphism of posets
with respect to the Bruhat order and the inclusion order
on ideals. 
\eco    
\subsection{}
\label{4.2}
We will need the following partial order on the set of subsets 
of $[1, m+n]$ with $k$ elements: for $J= \{ j_1 < \ldots < j_k \}$
and $J'= \{ j'_1 < \ldots < j'_k \} 
\subseteq [1,m+n]$ set 
\begin{equation}
\label{Jord}
J \leq J', \quad \mbox{if} \; \; j_l \leq j'_l \; \; \mbox{for all}
\; \; l=1, \ldots, k. 
\end{equation}
Set $J < J'$, if $J \leq J'$ and $J \neq J'$. 

For $J \subseteq [1, m+n]$ denote 
\[
p_1 (J) = J \cap [1,m] \; \; \mbox{and} \; \; 
p_2 (J) = J \cap [m+1, m+n].
\]
Given $J = \{ j_1, \ldots, j_k \} \subseteq [m+1, m+n]$, set 
$J-m = \{ j_1 - m, \ldots, j_k -m \} \subseteq [1,n]$.

For $J \subseteq [1, m+n]$, $|J|=k$ such that $J \leq c^m([1,k])$, define
\begin{equation}
\label{Delta}
\De^q(J) = \Delta^q_{w^\ci_m(p_1(J) \backslash p_1(c^m([1,k]))) , 
(p_2 (c^m([1, k])) \backslash p_2(J)) -m},
\end{equation}
where $w^\ci_m$ denotes the longest element of the copy of $S_m$ 
inside $S_{m+n}$ acting on the first $m$ indices.
First we simplify \eqref{Delta} in the cases $k \leq n$ and $k>n$, 
and verify that the two sets in the definition 
of the quantum determinant in the right hand side of \eqref{Delta}
have the same cardinality. Indeed, for $1 \leq k \leq n$ 
we have $p_1(c^m([1,k])) = \emptyset$ and $p_2(c^m([1,k])) = [m+1, m+k]$, 
and \eqref{Delta} simplifies to 
\begin{equation}
\label{Delta1} 
\De^q(J) =
\De^q_{w_m^\ci(p_1(J)), ([m+1, m+k] \backslash p_2(J)) -m}.
\end{equation}
Moreover for these values of $k$, 
$J \leq c^m([1,k]) =[m+1, m+k]$ implies 
$p_2(J) \subseteq [m+1, m+k]$ and thus
$|w_m^\ci(p_1(J))|=|([m+1, m+k] \backslash p_2(J)) -m|$.
Here and below $|S|$ denotes the cardinality of a set $S$.

For $n+1 \leq k \leq m+n$ we have 
$p_1(c^m([1,k])) = [1, k-n]$ and 
$p_2(c^m([1,k])) = [m+1, m+n]$, and \eqref{Delta} simplifies to 
\begin{equation}
\label{Delta2} 
\De^q(J)= 
\De^q_{w_m^\ci(p_1(J) \backslash [1,k-n]), 
([m+1, m+n] \backslash p_2(J)) -m}.
\end{equation}
For these values of $k$, 
$J \leq c^m([1,k]) = [1,k-n] \sqcup [m+1, m+n]$ 
implies $p_1(J) \supseteq [1,k-n]$, therefore
$|w_m^\ci(p_1(J) \backslash [1,k-n])|= 
|([m+1, m+n] \backslash p_2(J)) -m|$.

Denote by $\Lambda_q(\KK^{m+n})$ the quantum 
exterior algebra in $m+n$ generators.
It is a $\UU_q({\mathfrak{sl}}_{m+n})$-module algebra 
with generators $v_1, \ldots, v_{m+n}$
and relations $v_i v_j = - q v_j v_i$, 
$1 \leq j < i \leq m+n$, $v_i^2=0$, $i=1, \ldots, m+n$.
The $\UU_q({\mathfrak{sl}}_{m+n})$-action on it is given by:
\[
X_i^+ v_j = \delta_{i+1,j} v_i, 
\;
X_i^- v_j = \delta_{ij} v_{i+1},
\;
K_i v_j = q^{a_{ij}} v_j, 
\; \; i \in [1,m + n-1], j \in [1,m + n],
\]
where $a_{ij} = 1$ if $j=i$, 
$a_{ij}=-1$ if $j = i + 1$, and
$a_{ij}=0$ otherwise.  
The algebra $\Lambda_q(\KK^{m+n})$ is $\Zset$-graded by 
$\deg v_i = 1$. For $k =1, \ldots, m+n-1$ its 
component $\Lambda_q(\KK^{m+n})_k$ 
in degree $k$ is isomorphic to the fundamental 
representation $V(\om_k)$ of $\UU_q({\mathfrak{sl}}_{m+n})$. 
We will identify $V(\om_k)$ and $\Lambda_q(\KK^{m+n})$.
For $J=\{j_1 < \ldots < j_k\} \subseteq [1,m+n]$ define
\[
v_J = v_{j_1} \ldots v_{j_k}.
\]
When $J$ runs over all subsets of $[1,m+n]$ with $k$ elements 
we obtain a basis of $V(\om_k)$. 
The corresponding dual basis of 
$V(\om_k)^*$ will be denoted by $\{ \xi_J \mid J \subseteq [1, m+n], |J|=k \}$.
The Demazure modules $V_w(\om_k)$ are given by 
\[
V_w(\om_k) = \UU_+ T_w v_{[1,k]} = \Span
\{ v_J \mid J \subset [1,m+n], 
|J|=k, J \leq w([1,k]) \},
\]  
for all $w \in S_{m+n}$. For $J \subset [1,m+n]$, 
$|J|=k$, $J \leq w([1,k])$ denote
\[
\eta_J = \xi_J|_{V_{c^m}(\om_k)}.
\]
Then 
\begin{equation}
\label{dualDemazure}
\{ \eta_J \mid J \subset [1,m+n], 
|J|=k, J \leq c^m([1,k]) \}
\; \; \mbox{is a $\KK$-basis of} \; \; 
V_{c^m}(\om_k)^*.
\end{equation}
For a set $J$ as in \eqref{dualDemazure} and $y \in S_{m+n}^{\leq c^m}$
we have:
\begin{equation}
\label{inter}
\eta_J \in (V_w(\om_k)\cap \UU_- T_y v_{[1,k]})^\perp
\Leftrightarrow
J \ngeq y([1,k]).
\end{equation}  
Denote
\begin{multline}
\label{Upyw}
\Upsilon(y) = \{ J \subseteq [1,m+n] \mid |J|= k, k \in [1, m+n-1], \\
J \leq c^m([1,k]), J \ngeq y([1,k])  
\}.
\end{multline}

\ble{map} For all $k =1, \ldots, m+n$ and $J \subset [1, m+n]$
such that $|J|=k$ and $J \leq c^m([1,k])$ we have 
\begin{equation}
\label{coeff}
\zeta_{m,n}( \psi_{c^m}(d^{c^m,\om_k}_{\eta_J} ) ) = 
t' \De^q(J) 
\end{equation}
for some $t' \in \KK^*$ (depending on $J$).
\ele
\begin{proof} We will prove \leref{map} in the case $k \in [1,n]$. 
The case $k \in [n+1,m+n-1]$ is analogous and is left to the reader, cf. 
the proof of \cite[eq. (5.18)]{Y}.

Denote the longest element of $S_{m+n}$ by $w^\ci_{m+n}$ and consider 
the reduced decomposition 
\[
w^\ci_{m+n} = s_1 (s_2 s_1) \ldots (s_{m+n-1} \ldots s_1).
\]
Denote the Lusztig root vectors of $\UU^{w^\ci_{m+n}}_+= \UU_+$ 
from eq. \eqref{rootv} by 
\[
Y_{1,2}; Y_{1,3}, Y_{2,3}; \ldots; Y_{1,m+n}, \ldots, Y_{m+n-1, m+n}.
\]
By \cite[Lemma 2.1.1]{MC}, $Y_{i, i+1} = X_i^+$ for $i \in [1, m+n-1]$
and for $j>i+1$, $Y_{ij}$ is recursively given by
\begin{equation}
\label{rec1}
Y_{ij} = Y_{i, j-1}Y_{j-1,j} - q^{-1} Y_{j-1,j} Y_{i, j-1}.
\end{equation}
Using induction on $j-i$, one easily verifies that for all 
$I \subseteq [1, m+n]$, $|I|=k$ and $i, j \in [1,m+n]$, 
$i <j$:
\begin{align}
\label{ac}
(\tau Y_{ij}) v_I &= (-q)^{i-j+ | I \cap [i+1, j-1] | + 1 } 
v_{(I \backslash \{ j \}) \cup \{ i \} }, \; \;  
\mbox{if} \; \; j \in I \; \; \mbox{and}
\\ 
\label{ac1}
(\tau Y_{ij}) v_{I} & = 0,  \; \;  
\mbox{if} \; \; j \notin I,
\end{align}
recall \eqref{tau}.
The quantum $R$-matrix corresponding to 
$c^m \in S_{m+n}$ is given by 
\begin{align*}
\RR^{c^m} = &\left( \exp_q(q' X_{m, m+n} \otimes X_{m+n,m}) \dots 
\exp_q(q' X_{1, m+n} \otimes X_{m+n,1}) \right) \dots
\\
&\left( \exp_q(q' X_{m, m+1} \otimes X_{m+1,m}) \dots 
\exp_q(q' X_{1, m+1} \otimes X_{m+1,1}) \right),
\end{align*}
where $q' = q^{-1} - q$. Recall that $w^\ci_m$ denotes
the longest element of the copy of $S_m$ inside $S_{m+n}$ 
acting on the first $m$ indices. Fix $k \in [1,n]$ and $J \subset [1, m+n]$
such that $|J|=k$ and $J \leq c^m([1,k])$.
Denote
\[
w^\ci_m(p_1(J)) = \{ i_1 < \ldots < i_l \}, \; \; 
[m+1,m+n] \backslash p_2(J) = \{m + j_1 < \ldots < m + j_l \}.
\]
By \cite[Lemma 2.1.3]{MC}
$X_{i,m+j} = T^{-1}_{w^\ci_m} (Y_{i,m+j})$, $\forall i \in [1,m]$, 
$j \in [1,n]$. 
After some straightforward computations 
using eqs. \eqref{comp} and \eqref{ac}-\eqref{ac1} 
one deduces that 
\begin{align*}
&\zeta_{m,n}( \psi_{c^m}(d^{c^m,\om_k}_{\eta_J} ) )
= t'' \zeta_{m,n}( ( d^{c^m,\om_k}_{\eta_J} \tau \otimes \id) (\RR^{c^m}) ) 
\\
= & t'' \sum_{w \in S_l} 
\lcor \xi_{ \{i_1, \ldots, i_l\} }, 
(\tau Y_{i_{w(1)}, m+ j_1}) \ldots (\tau Y_{i_{w(l)}, m+ j_l}) 
v_{ \{m+ j_1, \ldots, m+ j_l \} } \rcor
x_{i_{w(l)}, j_l} \ldots x_{i_{w(1)}, j_1}
\\
= & t' \sum_{w \in S_l} (-q)^{- l(w)}
\lcor \xi_{ \{i_1, \ldots, i_l\} }, v_{ \{ i_1, \ldots, i_l \} } \rcor
x_{i_{w(l)}, j_l} \ldots x_{i_{w(1)}, j_1}
= t' \De^q(J)
\end{align*}
for some $t', t'' \in \KK^*$, which completes the proof of the lemma.
\end{proof}

Recall from \coref{RqMmn} that all $\Tset^{m+n}$-invariant prime ideals 
of $R_q[\Mmn]$ with respect to the action \eqref{torus2} are of the form
$I_{m,n}(y)$ for some $y \in S_{m+n}^{\leq c^m}$. The following theorem 
provides a constructive proof of the Goodearl--Lenagan conjecture \cite{GL},
showing that each of these ideals posses a polynormal generating sequence.
Even more, such a sequence is shown to be $P$-polynormal with respect to 
the action \eqref{newPact}, where $P$ is the weight lattice of $\sl_{m+n}$.

\bth{GL-thm} Assume that the base field $\KK$ has characteristic $0$ 
and $q \in \KK$ is transcendental over $\Qset$.
Fix any linear ordering $\prec$ on the set $\Upsilon(y)$
given by \eqref{Upyw} with the property that if $J, J' \in \Upsilon(y)$, 
$|J| = |J'|$ and $J < J'$ (recall \eqref{Jord}), 
then $J \prec J'$. Let 
$\Upsilon(y)= \{ J_1 \prec J_2 \prec \ldots \prec J_{\Upsilon(y)} \}$.
Then 
\begin{equation}
\label{genset}
\De^q(J_1), \ldots, \De^q(J_{|\Upsilon(y)|})
\end{equation}
is a $P$-polynormal generating sequence for the $\Tset^{m+n}$-invariant
prime ideal $I_{m,n}(y)$ of $R_q[\Mmn]$ (with respect to the action 
\eqref{newPact} of the weight lattice $P$ of $\sl_{m+n}$ on 
$R_q[\Mmn]$).
\eth
\begin{proof} Fix $y \in S_{m+n}$, $y \leq c^m$.
\thref{Yold} (c) and \leref{map} imply that the set from \eqref{genset} 
generates the ideal $I_{m,n}(y)$, see \cite[Theorem 5.5]{Y} for details.
By a direct computation one obtains
\[
\Wt(\eta_J) = - \om_k + \sum_{i=1}^{k} (\al_i + \ldots + \al_{j_i-1}), 
\]
for all $J=\{j_1 < \ldots < j_k \} \in \Upsilon(y)$.
In the setting of \thref{polynormal} we can choose
$\Om_{c^m}(y) = \{ \eta_J \mid J \in \Upsilon(y) \}$.  
The definition of the linear ordering $\prec$ implies that
\[
\eta_{J_1} < \ldots < \eta_{|\Upsilon(y)|}
\]
is a linear ordering on the set $\Om_{c^m}(y)$ satisfying
the conditions of \thref{polynormal}. 
It follows from \thref{polynormal} that 
\[
\psi_{c^m} \Big( d^{c^m, |J_1|}_{\eta_{J_1}} \Big), \ldots, 
\psi_{c^m} \Big( 
d^{c^m, |J_{|\Upsilon(y)|}|}_{\eta_{J_{|\Upsilon(y)|}}} \Big)
\]
is a $P$-polynormal generating sequence of the ideal $I_{c^m}(y)$
of $\UU_-^y$. Recall that $I_{m,n}(y) = \zeta_{m,n}(I_{c^m}(y))$,
see \eqref{ideal}. The statement of the theorem now follows from 
\eqref{coeff} and the fact that
$\zeta_{m,n} \colon \UU_-^{c^m} \to R_q[\Mmn]$
intertwines the $P$-actions \eqref{Pact} and \eqref{newPact}.
\end{proof}

\bre{GL_re} \thref{pol2} implies that under the weaker assumption 
that $q \in \KK^*$ is not a root of unity and without restrictions 
on the characteristic of the base field $\KK$, all $\Tset^{m+n}$-invariant 
prime ideals of $R_q[\Mmn]$ are $P$-polynormal with respect to 
the action \eqref{newPact}. (This is sufficient for our applications 
of polynormality, see Section \ref{catenarity} below.) 
Currently, there is no proof that under the above weaker 
assumptions the ideals $I_{m,n}(y)$ are generated by 
the quantum minors in \thref{GL-thm} and more generally 
that the ideals $I_w(y)$ of $\UU^w_-$  
are generated by the elements in \thref{polynormal}
(for an arbitrary simple Lie algebra $\g$). 
We conjecture that this is correct.
It is proved for the height one $\Tset^{m+n}$-primes 
of quantum matrices \cite[Proposition 4.2]{LLR} 
and more generally for the height one $\Tset^r$-prime ideals 
of all algebras $\UU^w_-$, \cite[Proposition 6.8]{Y4}.
\ere
\sectionnew{Catenarity of $\Spec \UU^w_-$} 
\label{catenarity}
\subsection{}
\label{5.1}
In this section we prove 
that $\Spec \UU^w_-$ is normally separated for all simple Lie
algebras $\g$ and $w \in W$.
From this we deduce that all algebras $\UU^w_-$ 
are catenary. Furthermore, we prove a formula 
for the heights of all $\Tset^r$-invariant prime ideals 
$I_w(y)$ of $\UU^w_-$, recall \thref{Yold}.

We recall that for a ring $R$, one says that $\Spec R$ is normally 
separated if for any two prime ideals $I \subsetneq I'$ of $R$ 
there exists $u \in I'$, which is normal in $R$ modulo $I$ 
and such that $u \notin I$. If, in addition a group $G$ acts on $R$ by algebra 
automorphisms, we say that $G-\Spec R$ is $G$-normally separated 
if for every two $G$-prime ideals $I \subsetneq I'$ of $R$ 
there exists $u \in I'$, which is $G$-normal in $R$ modulo $I$ 
and such that $u \notin I$, recall \S \ref{3.1}.

\bco{norm} Let $\KK$ be an arbitrary base field, $q \in \KK^*$ not a 
root of unity, $\g$ be a simple Lie algebra and $w \in W$.
Then $\Tset^r - \Spec \UU^w_-$ is $P$-normally separated
for the action \eqref{Pact}.
In particular, $\Tset^r - \Spec \UU^w_-$ is 
$\Tset^r$-normally separated with respect to the 
action \eqref{Tract}.
\eco

The special case of \coref{norm} for the algebras of quantum matrices
is due to Cauchon \cite{Cau}. The case when $w= w_0$ and $\KK= \Cset(q)$ 
is due to Caldero \cite{Cal}.
\\ \hfill \\
\noindent
{\em{Proof of \coref{norm}.}} Let $I \subsetneq I'$ be two $\Tset^r$-prime
ideals of $\UU^w_-$. By \thref{pol2}, $I'$ possesses 
a $P$-polynormal generating sequence $u_1, \ldots, u_n$.
Denote the image of $u_j$ in $\UU^w_-/I$ by $\ol{u}_j$,
$j = 1, \ldots, n$.
Let $k = \min \{ j \mid \ol{u}_j \neq 0 \}$. Then 
$u_k$ is $P$-normal modulo the ideal generated by 
$u_1, \ldots, u_{k-1} \in I'$. Therefore $u_k \in I'$ 
is a $P$-normal element of $\UU^w_-$ modulo $I$ and $u_k \notin I$. 
\qed
 
Next, we give a second proof to the $P$-normal separation 
result for the $\Tset^r$-invariant prime ideals 
$\{ I_w(y) \}_{y \in W^{\leq w} }$ of $\UU^w_-$, using
results of Gorelik \cite{G}. This proof also constructs 
explicit separating elements for all pairs of 
$\Tset^r$-invariant prime ideals of $\UU^w_-$. 

\bth{2pf} Assume that $\KK$ is an arbitrary base field, 
$q \in \KK^*$ is not a root of unity and $\g$ is an 
arbitrary simple Lie algebra. 
Let $y_1, y_2 \in W^{\leq w}$ and $\la' \in P^+$ 
be such that $y_1 < y_2$ in the Bruhat order and 
$y_1(\la') \neq y_2(\la')$. Then 
$\phi_w( e^{\la'}_{y_1} e^{-\la'}_w ) \in I_w(y_2)$ 
is a $P$-normal element of $\UU^w_-$ 
modulo $I_w(y_1)$, which does not belong to 
$I_w(y_1)$. For all
$\la \in P^+$, $\nu \in P$, $\xi \in (V(\la)^*)_\nu$, 
we have 
\begin{multline}
\label{c}
\phi_w( e^{\la'}_{y_1} e^{-\la'}_w ) 
\phi_w( c_\xi^\la e^{-\la}_w ) =
q^{ \lcor -(y_1 + w)(\la'), \nu +w(\la) \rcor }
\phi_w( c_\xi^\la e^{-\la}_w ) 
\phi_w( e^{\la'}_{y_1} e^{-\la'}_w )
\\
= \Big( ( - (y_1 + w)(\la')) .
\phi_w( c_\xi^\la e^{-\la}_w ) \Big) 
\phi_w( e^{\la'}_{y_1} e^{-\la'}_w )
\mod I_w(y_1).
\end{multline}
\eth
\begin{proof} \leref{comm-main} and the fact that 
$\phi_w \colon R^w_0 \to \UU^w_-$ is a graded 
antiisomorphism (see \thref{hom1}) 
imply \eqref{c}. Recall that 
\[
e^{\la'}_{y_1} = c^{\la'}_{\xi_{y_1, {\la'} }}
\]
for certain $\xi_{y_1, \la'} \in (V(\la')^*)_{- y_1 \la'}$, see 
\eqref{e}. Since $y_1 < y_2$ and $y_1(\la') \neq y_2 (\la')$, we 
have that $y_1(\la') > y_2 (\la')$. Therefore 
$(\UU_- T_{y_2} v_{\la'})_{y_1(\la')} = 0$ and 
$\xi_{y_1, {\la'} } \perp \UU_- T_{y_2} v_{\la'}$. 
Hence $e^{\la'}_{y_1} e^{-\la'}_w \in Q_w(y_2)^-$ and 
$\phi_w( e^{\la'}_{y_1} e^{-\la'}_w ) \in I_w(y_2)$.

Proposition 5.3.3 (ii), Lemma 6.6, and Corollary 6.10.1 (i)
of Gorelik \cite{G} imply that 
$\phi_w( e^{\la}_{y_1} e^{-\la}_w ) \notin I_w(y_1)$
for all $\la \in P^+$.  These results were formulated in 
\cite{G} for $\KK$ of characteristic $0$ and $q \in \KK$ 
transcendetal over $\Qset$, 
but it was shown in \cite[Theorem 3.1(b)]{Y5} that 
Gorelik's proof works for all fields $\KK$ and 
$q \in \KK^*$ which are not roots of unity.
This completes the proof of the theorem.
\end{proof}

We proceed with proving that all algebras $\UU^w_-$ have 
normal separation.

\bth{nsepar} Assume that $\KK$ is an arbitrary base field 
and $q \in \KK^*$ is not a root of unity. For all
simple Lie algebras $\g$ and $w \in W$, 
$\Spec \UU^w_-$ is normally separated.
\eth

Normal separation of $\Spec \UU^w_-$ was established in two special cases 
earlier. The case of the algebras of quantum matrices is 
due to Cauchon \cite{Cau}, who used very different techniques based 
on his method of deleting derivations. The case when $w = w_0$ 
(the longest element of the Weyl group $W$) and $\KK= \Cset(q)$
was obtained by Caldero in \cite{Cal}.
\\ \hfill \\
{\em{Proof of \thref{nsepar}.}} Goodearl proved \cite[Corollary 4.6]{G1} 
that, if $R$ is a right noetherian ring graded by an abelian 
group and $R$ has graded normal separation, then $\Spec R$ is 
normally separated. (A graded ring $R$ is said to have graded normal 
separation if for every two graded prime ideals $I \subsetneq I'$ 
there exists a homogeneous nonzero element $x \in I'/I$ 
which is normal in $R/I$.) The algebras $\UU^w_-$ are noetherian, 
because they are iterated skew polynomial rings. The 
graded prime ideals of $\UU^w_-$ with respect to the
$Q$-grading \eqref{Qgr} are precisely the $\Tset^r$-invariant 
prime ideals with respect to the action \eqref{Tract}.
\coref{norm} implies that the the set of $Q$-graded prime 
ideals of $\UU^w_-$ is $P$-normally separated. Recall that $P$-normal elements 
are $P$-eigenvectors (see \S \ref{3.1}) and that the $P$-eigenvectors 
in $\UU^w_-$ are precisely the homogeneous elements of $\UU^w_-$
with respect to the grading \eqref{Qgr}. 
Therefore the algebras $\UU^w_-$ have graded normal separation and
we can apply Goodearl's result to them, which establishes the theorem. 
\qed
\subsection{}
\label{5.2} We proceed with proving that all algebras $\UU^w_-$ 
are catenary. Motivated by Gabber's proof of catenarity of the 
universal enveloping algebras of all solvable Lie algebras, 
Goodearl and Lenagan proved the following theorem.

\bth{G-L} (Goodearl--Lenagan, \cite{GLen0}) Assume that $A$ is an 
affine, noetherian, Auslander--Gorenstein and Cohen--Macaulay
algebra over a field, with finite Gelfand--Kirillov dimension. 
If $\Spec A$ is normally separated, then $A$ is catenary. 
If, in addition, $A$ is a prime ring, then Tauvel's height 
formula holds.
\eth
We recall that Tauvel's height formula holds for $A$ 
if for all prime ideals $I$ of $A$, the height of $I$ 
is equal to 
\[
\GKdim A - \GKdim(A/I).
\]
Here and below $\GKdim(.)$ denotes the Gelfand--Kirillov dimension 
of an algebra or a module.

For the convenience of the reader we also recall the notions 
of Auslander regular, Auslander--Gorenstein, and Cohen--Macaulay
rings. A ring $R$ is called Auslander--Gorenstein if 
the injective dimension of $R$ (as both right and left 
$R$-module) is finite, and for all integers $0 \leq i < j$ 
and finitely generated (right or left) $R$-modules $M$, 
we have $\Ext^i_R(N, R) = 0$ 
for all $R$-submodules $N$ of $\Ext^j_R(M,R)$. A ring $R$ is 
said to be Auslander regular if, in addition, the global dimension 
of $R$ is finite. The grade of a finitely generated
$R$-module $M$ is defined by
\[
j(M) = \inf \{i \geq 0 \mid \Ext^i_R(M,R) \neq 0 \}.
\]
An algebra $R$ is called Cohen--Macauley if 
\[
j(M)+ \GKdim M = \GKdim R 
\]
for all finitely generated $R$-modules $M$. 

We apply \thref{G-L} to the algebras $\UU^w_-$. The normal separation 
of $\Spec \UU^w_-$ was established in \thref{nsepar}. It is well 
known that all algebras $\UU^w_-$ are CGL extensions 
(a special kind of iterated skew polynomial 
rings) in $l(w)$ variables, where $l(w)$ is the length 
of $w$. Thus for all base fields $\KK^*$, $q \in \KK^*$ 
not a root of unity, and $w \in W$, the algebra 
$\UU^w_-$ are affine, noetherian and 
\begin{equation}
\label{GKU}
\GKdim \UU^w_- = l(w).
\end{equation}

We derive from the following result of Ekstr\"om, 
Levasseur and Stafford 
\cite{E,LS} that all algebras $\UU^w_-$ are 
Auslander regular and Cohen--Macaulay:

\bpr{LS} (Ekstr\"om, Levasseur--Stafford, \cite{E,LS}) Assume $R$ is a 
noetherian, Auslander regular ring. Let $S= R[x; \sigma, \delta]$ 
be a skew polynomial extension of $R$. Then: 

(a) $S$ is Auslander regular.

(b) If $R = \oplus_{k \geq 0} R_k$ is a connected graded 
Cohen--Macauley $\KK$-algebra over a field $\KK$ such 
that $\sigma (R_k) \subseteq R_k$ 
for all $k \geq 0$, then $S$ is Cohen--Macauley. 
\epr 
\bpr{Uwhom} For all base fields $\KK$, $q \in \KK^*$ 
not a root of unity, simple Lie algebras $\g$, and Weyl group 
elements $w \in W$, the algebras $\UU^w_-$ are 
Auslander regular and Cohen--Macauley.
\epr
\begin{proof} Fix $w \in W$ and a reduced expression
$w= s_{i_1} \ldots s_{i_l}$ of it. Denote $w' = w s_{i_l}$. 
In the notation of \S 2.2, the subalgebra of $\UU^w_-$ 
generated by $X^-_{\be_1}, \ldots, X^-_{\be_{l-1}}$ 
coincides with $\UU^{w'}_-$. Recall the 
Levendorskii--Soibelman straightening rule
\begin{multline}
\label{LS}
X^-_{\be_i} X^-_{\be_j} - q^{\lcor \be_i, \be_j \rcor }
X^-_{\be_j} X^-_{\be_i}  \\
= \sum_{ {\bf{n}} = (n_{i+1}, \ldots, n_{j-1}) \in \Nset^{\times (j-i-2)} }
p_{\bf{n}} (X^-_{\be_{i+1}})^{n_{i+1}} \ldots (X^-_{\be_{j-1}})^{n_{j-1}},
\; \; p_{\bf{n}} \in \KK,
\end{multline}
for all $1 \leq i < j \leq l$. It implies that $\UU^w_-$ 
is a skew polynomial extension of $\UU^{w'}_-$:
\[
\UU^w_- \cong \UU^{w'}_- [x; \sigma, \delta], 
\]
where $\sigma \in \Aut ( \UU^{w'}_-)$ is given by 
\begin{equation}
\label{sigma}
\sigma(x) = q^{- \lcor \ga , \be_l \rcor } x, \quad 
\forall x \in (\UU^{w'}_-)_{- \ga}, \ga \in Q^+,
\end{equation}
recall \eqref{Tract}. By repeated applications of \prref{LS} (a)
one obtains that the algebra $\UU^w_-$ is Auslander regular.
Fix $\mu \in \sum_{i=1}^r \Zset_+ \om_i$ and 
specialize the $-Q^+$-grading of $\UU^{w'}_-$ to an 
$\Nset$-grading by
\[
(\UU^{w'}_-)_k = \oplus_{\mu \in Q^+} 
\{ (\UU^{w'}_-)_{- \ga}  \mid \lcor \ga , \mu \rcor = k \}, 
\quad
k \in \Nset.
\]
Obviously $\UU^w_-$ is connected. It follows 
from \eqref{sigma} that 
$\sigma((\UU^{w'}_-)_{- \ga} ) \subseteq (\UU^{w'}_-)_{- \ga}$
for all $\ga \in Q^+$ and thus 
$\sigma((\UU^{w'}_-)_k) \subseteq (\UU^{w'}_-)_k$
for all $k \in \Nset$. Repeatedly applying \prref{LS} (b),  
we obtain that the algebra $\UU^w_-$ is Cohen--Macauley.
\end{proof}

Combining Theorems \ref{tnsepar} and \ref{tG-L}, \prref{Uwhom},
and the fact that the algebras $\UU^w_-$ are affine, 
noetherian domains of Gelfand--Kirillov dimension $l(w)$,  
we obtain:
\bth{catenary} For all base fields $\KK$, $q \in \KK^*$ 
not a root of unity, simple Lie algebras $\g$ and Weyl group 
elements $w \in W$, the algebras $\UU^w_-$ are catenary
and Tauvel's height formula holds.
\eth

The special case of \thref{catenary} for the algebras of quantum matrices
is due to Cauchon \cite{Cau}. The case when $w=w_0$ and $\Char \KK =0$ 
was obtained by Malliavin \cite{M}, Goodearl and Lenagan \cite{GLen0}.
\subsection{}
\label{5.3}
In this subsection we establish formulas for the heights 
of all $\Tset^r$-invariant prime ideals $I_w(y)$ of $\UU^w_-$ 
(recall \thref{Yold}) and the Gelfand--Kirillov dimensions 
of the quotients $\UU^w_-/I_w(y)$.

\bth{height} For all base fields $\KK$, $q \in \KK^*$ 
not a root of unity, simple Lie algebras $\g$ and Weyl group 
elements $y, w \in W$, $y \leq w$, the height of the 
$\Tset^r$-invariant prime ideal $I_w(y)$ equals $l(y)$ and
\[
\GKdim (\UU^w_-/I_w(y)) = l(w) - l(y).
\]
\eth

We will need the following proposition. Its proof was communicated to 
us by Ken Goodearl \cite{G2}.
\bpr{he} Let $A$ be a noetherian algebra over an infinite
field $\KK$, equipped with a rational action of a $\KK$-torus $T$ 
by algebra automorphisms. If $T-\Spec A$ is $T$-normally separated,
then for each pair of $T$-invariant prime ideals $I \subsetneq I'$ 
there exists a saturated chain of prime ideals 
$I \subsetneq I_1 \subsetneq \ldots \subsetneq I_m \subsetneq I'$
consisting entirely of $T$-invariant prime ideals.
\epr
\begin{proof} Arguing by induction, it suffices to prove that:

{\em{If $I \subsetneq I'$ is a pair of $T$-invariant prime ideals 
of $A$ such that there is no $T$-invariant prime ideal
of $A$ lying strictly between $I$ and $I'$, then there is no 
prime ideal of $A$ lying strictly between $I$ and $I'$.}}

By changing $A$ to $A/I$, we see that it is sufficient to prove the 
special case when $I= \{0\}$ (and thus $A$ is a prime algebra).
The assumption on normal separation 
implies that there exists a normal element $c$ of $A$ such that
$c \in I'$, $c \neq 0$ and $c$ is a $T$-eigenvector. 
By \cite[Proposition II.2.9]{BG} 
all $T$-primes of $A$ are prime. Therefore all minimal primes over 
$c A$ are $T$-invariant prime ideals. Denote one of them 
that is contained in $I'$ by $J$. 
We have that $\{0\} \subsetneq J \subseteq I'$. The assumption that
there is no $T$-invariant prime ideal of $A$ lying strictly 
between $\{0\}$ and $I'$ implies that $J= I'$. By the 
principle ideal theorem 
\cite[Theorem 4.1.11]{McR} the height of $J$ equals $0$ 
or $1$. The former is impossible since $\{0\} \subsetneq I'$
and $\{0\}$ is a prime ideal. 
Thus the height of $I'$ is equal to $1$, and there are 
no prime ideals of $A$ lying strictly between $\{0\}$ and 
$I'$.
\end{proof}
\noindent
{\em{Proof of \thref{height}.}} By \coref{norm}, $\Tset^r - \Spec \UU^w_-$ 
is $\Tset^r$-normally separated. The base field $\KK$ is infinite since
$q \in \KK^*$ is not a root of unity. Applying \prref{he}, we obtain 
that there exists a saturated chain of prime ideals
$ \{ 0 \} \subsetneq I_1 \subsetneq \ldots \subset I_m \subsetneq I_y(w)$
consisting entirely of $\Tset^r$-invariant prime ideals. It follows from 
\thref{Yold} (b) that the length of this chain is equal to $l(y)$. Therefore
the height of $I_y(w)$ equals $l(y)$. Recall from \eqref{GKU} that 
$\GKdim \UU^w_- = l(w)$. Applying the fact 
that Tauvel's height formula holds for $\UU^w_-$ (\thref{catenary}), 
we obtain that
\[
\GKdim ( \UU^w_- /I_w(y) ) = \GKdim \UU^w_- - l(y) = l(w) - l(y),
\]      
which completes the proof of the theorem.
\qed

\end{document}